\documentclass[a4paper,11pt,oneside]{article}

\usepackage{amsmath,amsthm,amsfonts,amssymb,amscd}
\usepackage{mathtools}
\usepackage{float}

\usepackage{setspace}
\usepackage{enumerate}
\usepackage{booktabs}
\usepackage[noend]{algpseudocode}
\usepackage[normalem]{ulem}
\usepackage{graphicx,bm,xcolor}
\usepackage{algorithm}
\usepackage{algpseudocode}
\usepackage{hyperref} 
\usepackage{caption}
\usepackage[T1]{fontenc}
\usepackage{t1enc}
\usepackage[polish,german,english]{babel}
\usepackage{verbatim}
\usepackage{diagbox}
\usepackage{colortbl}
\usepackage{siunitx}
\usepackage{mathdots} 
\usepackage{bbm}
\usepackage{cancel} 
\usepackage{subcaption}
\usepackage{pgf,tikz,pgfplots}
\usepackage{mathrsfs}
\usetikzlibrary{arrows}
\usetikzlibrary{shapes.misc}
\usetikzlibrary{shapes.geometric, arrows,chains}
\usetikzlibrary{positioning}
\usetikzlibrary{decorations.pathreplacing} 
\usetikzlibrary{calc}
\usepackage{tikz-cd} 

\usepackage{xcolor} 

\usepackage{url}
\usepackage{authblk}

\algnewcommand{\algorithmicgoto}{\textbf{go to}}
\algnewcommand{\Goto}[1]{\algorithmicgoto~\ref{#1}}

\newcommand{\R}{\mathbb{R}}

\definecolor{orange}{RGB}{230, 159, 0}
\definecolor{skyblue}{RGB}{86, 180, 233}
\definecolor{yellow}{RGB}{240, 228, 66}
\definecolor{blue}{RGB}{0, 114, 178}
\definecolor{vermillion}{RGB}{213, 94, 0}

\usepackage[top=2cm, bottom=2cm, left=2cm, right=2cm]{geometry}

\theoremstyle{plain} 
\newtheorem{theorem}{Theorem}[section]
\newtheorem{proposition}{Proposition}[section]

\newtheorem{remark}[theorem]{Remark}

\theoremstyle{definition} %

\theoremstyle{remark} %

\DeclareMathOperator*{\argmax}{arg\,max}

\DeclareMathOperator*{\dG}{dG}

\DeclareMathOperator*{\Ieff}{I_{\text{eff}}}
\DeclareMathOperator*{\Iind}{I_{\text{ind}}}

\newcommand\restrict[1]{\raisebox{-.5ex}{$|$}_{#1}}

\definecolor{brewer1}{HTML}{A6CEE3}
\definecolor{brewer2}{HTML}{1F78B4}
\definecolor{brewer3}{HTML}{B2DF8A}
\definecolor{brewer4}{HTML}{33A02C}

\begin{document}


\title{
Adaptive space-time model order reduction with dual-weighted residual (MORe DWR) error control for poroelasticity 
}

\author[1,2]{Hendrik Fischer}
\author[1,2]{Julian Roth}
\author[2]{Ludovic Chamoin}
\author[2]{Amélie Fau}
\author[3]{Mary Wheeler}
\author[1,2]{Thomas Wick}

\affil[1]{Leibniz Universit\"at Hannover, Institut f\"ur Angewandte
  Mathematik, \linebreak AG Wissenschaftliches Rechnen, Welfengarten 1, 30167 Hannover, Germany}
\affil[2]{Universit\'e Paris-Saclay, CentraleSupélec, ENS Paris-Saclay, CNRS, LMPS - Laboratoire de Mécanique Paris-Saclay,
91190 Gif-sur-Yvette, France}
\affil[3]{The University of Texas at Austin, Oden Institute, Austin, TX 78712, USA}

\date{}
\maketitle


\begin{abstract}
In this work, the space-time MORe DWR (\underline{M}odel \underline{O}rder \underline{Re}duction with \underline{D}ual-\underline{W}eighted \underline{R}esidual error estimates) 
framework
is extended and further developed for single-phase flow problems in porous media. 
Specifically, our problem statement is the Biot system which consists of 
vector-valued displacements (geomechanics) coupled to a Darcy flow pressure 
equation.
The MORe DWR method introduces a goal-oriented adaptive incremental proper orthogonal decomposition (POD) based-reduced-order model (ROM). 
The error in the reduced goal functional is estimated during the simulation, and the POD basis is enriched on-the-fly if the estimate exceeds a given threshold. 
This results in a reduction of the total number of full-order-model solves for the simulation of the porous medium, a robust estimation of the quantity of interest and well-suited reduced bases for the problem at hand.
We apply a space-time Galerkin discretization with Taylor-Hood elements in space and 
a discontinuous Galerkin method with piecewise constant functions in time. The latter is 
well-known to be similar to the backward Euler scheme. 
We demonstrate the efficiency of our method on 
the well-known two-dimensional Mandel benchmark and a three-dimensional footing problem. 
\end{abstract}


\section{Introduction}
\label{sec_intro}
Porous media problems have long-standing applications in subsurface modeling, 
groundwater flow, hydraulic fracturing, geothermal energy recovery, and nuclear 
waste storage. 
The resulting mathematical models yield coupled systems 
of partial differential equations (PDEs) from which one of the most well-known is the 
so-called Biot problem \cite{Biot41a,Biot41b,Biot55,Biot7172,To92}. Developments over the last two decades have led 
to more complicated models, due to additional physics that can be incorporated. Nonstationary, nonlinear coupled PDE systems are often obtained. They are also sometimes subject 
to inequality constraints as in multiphysics phase-field fracture in porous media \cite{WheWiLee20,Wick2020PFF}. Despite advances in numerical solvers of iterative or multigrid type, high-performance 
parallel computing, and improved hardware performances, the computational cost 
for solving such multiphysics problems remains high. 
Specific examples in the field of poroelasticity are 
\cite{Bause2017,KimTchJua11b,MiWhe12,WHITE201655,WhiBo11,BASTIAN1999199,Molenar95,HoKrLyPh19,anselmann2023energyefficient,BOTH2017101}.

Model order reduction (MOR) and reduced-order modeling (ROM) techniques \cite{lassila_model_2014,benner2020model,sirovich1987turbulence,willcox2002balanced,Ekre2020,BeCoOhlWill15, Peherstorfer2015, Peherstorfer2016} 
provide one possible solution to significantly reduce the computational cost. Therein, the problem is split into two phases, an offline phase (solving a costly original high-fidelity model) to build a representative reduced basis and an online phase in which the reduced-order model 
is solved fast. This splitting comes with an additional cost, which pays 
off when the original model needs to be solved several (hundreds and more) times as 
for example in uncertainty quantification, inverse modeling, and optimal control. However, providing a robust and problem-specific reduced basis for nonlinear and general problems is a challenge.

Other model-order reduction approaches are based on solving the Galerkin problem with a separation of variables \cite{pgd_optimality}. For nonlinear problems, some robust methods have been proposed, but they require a highly intrusive and specific numerical framework \cite{latin_int}, or less intrusive and more flexible approaches, but without a mathematical proof of convergence for nonlinear cases \cite{latin_non_int}. Thus, there is a need for reduced-order numerical strategies, which are flexible, general and for which the error accompanying the reduced solution space can be controlled for any quantity of interest.

Thanks to our own recent advances in space-time modeling and goal-oriented 
space-time error control for nonstationary, coupled problems \cite{RoThiKoeWi2023}, 
an incremental error-controlled proper orthogonal decomposition (POD) based ROM method, i.e. 
MORe DWR (\underline{M}odel \underline{O}rder \underline{Re}duction with \underline{D}ual-\underline{W}eighted \underline{R}esidual error estimates), was 
suggested for the heat and elastodynamics equations in \cite{FiRoWiChaFau2023}. 
Therein, satisfying results for adaptively refining the ROM basis according to 
some distributed-in-time goal functionals were obtained.

In the current work, the Biot system of poroelasticity, namely vector-valued displacements and a Darcy-type pressure equation, is considered and modeled 
in a space-time fashion. Specifically, the temporal discretization is based 
on discontinuous Galerkin (dG) finite elements, while the spatial discretization 
employs classical continuous Galerkin (cG) finite elements. This serves 
as a basis for space-time error-controlled adaptivity by applying the dual-weighted 
residual method (DWR) \cite{becker_rannacher_2001,bangerth_rannacher_2003} 
to our ROM concepts. Importantly, the ROM error is computed with 
respect to a goal functional, which is often motivated from technical quantities 
of interest arising in physics, engineering and practical applications. 
The ROM updates 
are performed in an incremental fashion as recently suggested 
in \cite{kuhl2023incremental} and do not 
require a full FOM (full-order model) run in advance. Thus, the ROM is updated on the fly using FOM solutions estimated only for one temporal element. 
Here, the singular value decomposition (SVD) for computing the leading components
reduces to a truncated version 
and provides the snapshot matrix.
An important point in coupled problems resorts to considering a monolithic version 
and computing one single SVD or, alternatively, 
splitting the SVD into its PDE components. In this work, 
we aim for the latter and split the displacements and pressure into separate 
SVDs in order to facilitate computing different POD sizes for the individual 
solution components. These developments result into a final MORe DWR algorithm,
which is tested for two-dimensional and three-dimensional configurations. First,
the well-known Mandel benchmark is addressed. In a second numerical test a 
poroelastic three-dimensional block is considered. Therein, on the one hand, typical 
cost factors such as speedup, FOM solves and ROM sizes are investigated. 
On the other hand, standard quantities of a posteriori error control such as 
effectivity indices and indicator indices are computed as well, demonstrating 
the performance of our overall proposed concepts.

The outline of this paper is as follows. In Section~\ref{sec_problem_formulation}
the strong form problem statement is presented first. Then, a space-time formulation 
is derived. The dual-weighted residual method for goal-oriented error control  
is exposed, which
requires a linear, backward-in-time running, adjoint problem.
Next, in Section~\ref{sec_MOReDWR}, the MORe DWR (\underline{M}odel \underline{O}rder \underline{Re}duction with \underline{D}ual-\underline{W}eighted \underline{R}esidual error estimates) concepts are explained. 
Goal-oriented 
error control for refining the ROM basis is also discussed. In Section~\ref{sec:numerical_tests} two numerical examples are investigated. The first test is the two-dimensional 
Mandel benchmark. The second test is a three-dimensional footing problem.

\section{Problem formulation}
\label{sec_problem_formulation}
In this work, we model a porous medium by coupling a vector-valued 
displacement equation (geomechanics) to Darcy flow in a poro-elastic medium.
In the problem description, $I := (0,T)$ denotes the temporal
domain and $\Omega \subset \mathbb{R}^d$
a sufficiently smooth spatial domain with spatial dimension $d \in \{2,3\}$.

\subsection{Strong form of poroelasticity}
The governing equations for poroelasticity read \cite{Coussy2004,LeSchref99}: Find pressure $p: \bar{\Omega} \times \bar{I} \rightarrow \mathbb{R}$ and displacement $u: \bar{\Omega} \times \bar{I} \rightarrow \mathbb{R}^d$ such that
\begin{equation}\label{eq:strong_form_poroelasticity}
\begin{aligned}
    \partial_t(cp + \alpha(\nabla_x \cdot u)) - \frac{1}{\nu}\nabla_x \cdot (K\nabla_x p) &= 0 \qquad\quad \text{in } \Omega \times I, \\
    -\nabla_x \cdot \sigma(u) + \alpha \nabla_x p &= 0 \qquad\quad \text{in } \Omega \times I,
\end{aligned}
\end{equation}
with the isotropic stress tensor
\begin{align*}
    \sigma(u) := \mu(\nabla_x u + (\nabla_x u)^T) + \lambda (\nabla_x \cdot u)I,
\end{align*}
and the (constrained specific) storage coefficient $c \geq  c^* > 0$, 
which is assumed strictly positive in this paper,
the Biot-Willis constant $\alpha \in [0,1]$, the permeability tensor $K$, fluid’s viscosity $\nu$ and the Lamé parameters $\lambda, \mu > 0$.
We notice that the constant $c$ may depend on space, i.e., $c(x)$, and is linked to the compressibility 
$M>0$. When $c$ tends to zero, numerical instabilities may arise (e.g., \cite{Phi05,Phillips2009}) for which mixed methods or enriched Galerkin 
or discontinuous Galerkin finite elements in the numerical discretization should be employed \cite{Liu04,SunLiu09,MuThoLou96,Phillips2008,LeeLeeWhe16}.
This coupled system of equations is also known as the Biot system \cite{Biot41a,Biot41b,Biot55,Biot7172,To92}. 
A rigorous mathematical analysis of this poroelasticity problem can be found in \cite{Showalter2000}.

Assuming homogeneous Dirichlet boundary conditions for the displacement on $\Gamma_D$ and inhomogeneous Neumann/traction boundary conditions 
\begin{align*}
    \sigma(u) \cdot n  &= t
\end{align*}
on $\Gamma_N = \partial\Omega\setminus \Gamma_D$, the spatial function spaces are given by
\begin{align*}
   V_u(\Omega) &= \left(H^1_{0, \Gamma_D}(\Omega)\right)^d, \\
   V_p(\Omega) &= H^1(\Omega), \\
   V(\Omega) &= V_u(\Omega) \times V_p(\Omega).
\end{align*}
All boundary conditions, including pressure conditions as well, are carefully specified in 
Section \ref{sec:numerical_tests}.

\subsection{Weak space-time form of the primal problem}
We define the space-time function space $X(I, V(\Omega))$ as
\begin{align*}
    X(I, V(\Omega)) := L^2(I, V(\Omega)) \cap H^1(I, V^\ast(\Omega))
\end{align*}
with the dual space of $V(\Omega)$ being denoted as $V^\ast(\Omega) = L(V(\Omega), \mathbb{R})$.
Furthermore, we will use the notation
\begin{align*}
    (f,g) := (f,g)_{L^2(\Omega)} := \int_\Omega f \cdot g\ \mathrm{d}x, \qquad (\!(f,g)\!) := (f,g)_{L^2(I, L^2(\Omega))} := \int_I (f, g)\ \mathrm{d}t, \\
    \langle f,g \rangle := \langle f,g \rangle_{L^2(\Gamma)} := \int_\Gamma f \cdot g\ \mathrm{d}s, \qquad (\!\langle f,g\rangle\!) := (f,g)_{L^2(I, L^2(\Gamma))} := \int_I \langle f, g\rangle\ \mathrm{d}t.
\end{align*}
In this notation, $f \cdot g$ represents the Euclidean inner product if $f$ and $g$ are scalar- or vector-valued and represents the Frobenius inner product if $f$ and $g$ are matrices.

We can now derive the space-time variational formulation for the (primal) poroelasticity problem. By integration by parts, the variational formulation reads: Find $U := \{u, p\} \in X\left(I, V(\Omega)\right)$ such that
\begin{align*}
    c(\!(\partial_t p, \phi^p)\!) + \alpha (\!(\partial_t (\nabla_x \cdot u), \phi^p)\!) + \frac{K}{\nu}(\!(\nabla_x p, \nabla_x\phi^p)\!) &= 0, \\
    (\!(\sigma(u), \nabla_x\phi^u)\!) - \alpha(\!(pI, \nabla_x\phi^u)\!) + \alpha (\!\langle pn, \phi^u\rangle\!)_{ \Gamma_N \times I} &= (\!\langle t, \phi^u \rangle\!)_{ \Gamma_N \times I}\\
    \forall \Phi := \begin{pmatrix}
    \phi^u \\ 
    \phi^p 
    \end{pmatrix}  \in X\left(I, V(\Omega)\right).
\end{align*}
Here, $n$ denotes the normal vector on the Neumann boundary $\Gamma_N$.

For the time discretization, let $\mathcal{T}_k := \{I_m \mid 1 \leq m \leq M\}$ with $I_m := (t_{m-1}, t_m)$ and 
\begin{align*}
    0 =: t_0 < t_1 < \cdots < t_{M-1} < t_M =: T
\end{align*}
be a partitioning of time i.e., $\bar{I} = [0,T] = \cup_{m = 1}^M \bar{I}_m$.
Then, in an intermediate step similar to spatial discontinuous Galerkin ($\dG$) method, we can define a broken function space
that allows for discontinuities at the temporal grid points:
\begin{align*}
    \widetilde{X}(\mathcal{T}_k,V(\Omega)) := \left\lbrace
    \begin{pmatrix} u \\ p \end{pmatrix} \in L^2(I,L^2(\Omega)^{d+1}), 
    \begin{pmatrix} u \\ p \end{pmatrix}\middle|_{I_m} \in X(I_m,V(\Omega))\;\forall I_m\in\mathcal{T}_k
    \right\rbrace.
\end{align*}
Due to these discontinuities, the limits of a function $f$ at time $t_m$ from above and from below are
\begin{align*}
    f_{m}^\pm := \lim_{\epsilon \searrow 0} f(t_m \pm \epsilon),
\end{align*}
and the jump of the function value of $f$ at time $t_m$ is
\begin{align*}
    [f]_m := f_{m}^+ - f_{m}^-.
\end{align*}

\noindent Accounting for the discontinuities in time, we thus need to solve the problem: 

Find $U \in \tilde{X}\left(\mathcal{T}_k, V(\Omega)\right)$ such that
\begin{align}\label{eq:space_time_form_poroelasticity}
    A(U)(\Phi)  = F(\Phi) \qquad \forall \Phi \in \tilde{X}\left(\mathcal{T}_k, V(\Omega)\right).
\end{align}
The bilinear form and right-hand side read
\begin{align*}
    A(U)(\Phi) &= \sum_{m = 1}^M \int_{I_m}  c(\partial_t p, \phi^p) + \alpha (\partial_t (\nabla_x \cdot u), \phi^p) + \frac{K}{\nu}(\nabla_x p, \nabla_x\phi^p) \ \mathrm{d}t \\
    &+ \sum_{m = 1}^M \int_{I_m} (\sigma(u), \nabla_x\phi^u) - \alpha(pI, \nabla_x\phi^u) + \alpha\langle pn, \phi^u\rangle_{ \Gamma_N} \ \mathrm{d}t \\
    &+ \sum_{m = 1}^{M-1}\alpha([\nabla_x \cdot u]_{m}, \phi_{m}^{p,+}) + \alpha (\nabla_x \cdot u_{0}^{+},\phi_0^{p,+})
    + \sum_{m = 1}^{M-1}c([p]_{m}, \phi_{m}^{p,+}) + c (p_{0}^{+},\phi_0^{p,+})
\end{align*}
and 
\begin{align*}
    F(\Phi) &:= (\!\langle t, \phi^u\rangle\!)_{ \Gamma_N \times I} + \alpha(\nabla_x \cdot u^0, \phi_0^{p,+}) + c(p^0, \phi_0^{p,+}),
\end{align*}
where
\begin{align*}
    U := (u, p), \quad \Phi := (\phi^u, \phi^p).
\end{align*}
\begin{remark}
    This space-time formulation can also be easily extended to a multirate-in-time setting with different timestep sizes for displacement and pressure \cite{roth2023monolithic}; 
for a rigorous analysis starting directly with the backward Euler time discretization and mixed spaces for flow and conformal Galerkin for geomechanics is provided in \cite{Almani2016}.
\end{remark}

For the projection of the problem onto the POD vectors, the block structure of the linear system is exploited to build separate bases for displacement and pressure.
We can decompose the bilinear form into the four following blocks
\begin{align*}
    A_1(u)(\phi^u) &= \sum_{m = 1}^M \int_{I_m} (\sigma(u), \nabla_x\phi^u)\ \mathrm{d}t, \\
    A_2(p)(\phi^p) &= \sum_{m = 1}^M \int_{I_m}  c(\partial_t p, \phi^p) + \frac{K}{\nu}(\nabla_x p, \nabla_x\phi^p) \ \mathrm{d}t + \sum_{m = 1}^{M-1}c([p]_{m}, \phi_{m}^{p,+}) + c (p_{0}^{+},\phi_0^{p,+}), \\
    B_1(p)(\phi^u) &= \sum_{m = 1}^M \int_{I_m} - \alpha(pI, \nabla_x\phi^u) + \alpha\langle pn, \phi^u\rangle_{\Gamma_N} \ \mathrm{d}t, \\
    B_2(u)(\phi^p) &= \sum_{m = 1}^M \int_{I_m} \alpha (\partial_t (\nabla_x \cdot u), \phi^p)\ \mathrm{d}t + \sum_{m = 1}^{M-1}\alpha([\nabla_x \cdot u]_{m}, \phi_{m}^{p,+}) + \alpha (\nabla_x \cdot u_{0}^{+},\phi_0^{p,+}).
\end{align*}

\subsection{dG(0), i.e., backward Euler, time discretization of the primal problem}\label{sec:backward_euler_primal}
Using the space-time formulation of the poroelasticity problem (\ref{eq:space_time_form_poroelasticity}), a $\dG(0)$ time discretization is derived by using discontinuous, piecewise-constant finite elements in time; see e.g., \cite{Wi23_st}. This results into a backward Euler scheme.
Then, on each temporal element $I_m = (t_{m-1}, t_m)$, we have temporally constant functions $(u_m, p_m) \in V(\Omega)$ such that
\begin{align*}
    u\restrict{I_m} =: u_m, \qquad p\restrict{I_m} =: p_m. 
\end{align*}
Using piecewise-constant functions $f$ in time,
we can insert the relations
\begin{align*}
    \partial_t f\restrict{I_m} = 0, \qquad [f]_m = f_{m+1} - f_m
\end{align*}
and arrive at the time-stepping scheme: Find $U_m := (u_m, p_m) \in V(\Omega)$ such that
\begin{align*}
    &c(p_m-p_{m-1}, \phi^p) + \alpha(\nabla_x \cdot u_m-\nabla_x \cdot u_{m-1}, \phi^p) + k\frac{K}{\nu}(\nabla_x p_m, \nabla_x \phi^p) \\
    &+ k(\sigma(u_m), \nabla_x\phi^u) - k\alpha(p_mI, \nabla_x\phi^u) + k\alpha\langle p_m n, \phi^u\rangle_{\Gamma_N} = k\langle t, \phi^u\rangle_{\Gamma_N} \qquad \forall \Phi = (\phi^u, \phi^p) \in V(\Omega).
\end{align*}
Here, $k := t_m - t_{m-1}$ denotes a (constant) timestep size.

\begin{remark}\label{remark:backward_euler_primal}
    Instead of starting from the space-time formulation, we can also directly derive the backward Euler time discretization from the strong formulation of poroelasticity (\ref{eq:strong_form_poroelasticity}). The main difference is that now we do not scale 
    the mechanics equation by the timestep size $k$. This is motivated by the displacement equation being quasi-static. Hence, in traditional time discretizations the time integral is omitted for this equation. For a consistent mathematical description, the previous time-stepping scheme fits better in the space-time setting, but in the actual implementation we can neglect the time step size $k$ in the mechanics equation.  We then have:
    Find $U_m := (u_m, p_m) \in V(\Omega)$ such that
    \begin{align*}
        &c(p_m-p_{m-1}, \phi^p) + \alpha(\nabla_x \cdot u_m-\nabla_x \cdot u_{m-1}, \phi^p) + k\frac{K}{\nu}(\nabla_x p_m, \nabla_x \phi^p) \\
        &+ (\sigma(u_m), \nabla_x\phi^u) - \alpha(p_mI, \nabla_x\phi^u) + \alpha\langle p_m n, \phi^u\rangle_{\Gamma_N} = \langle t, \phi^u\rangle_{\Gamma_N} \qquad \forall \Phi = (\phi^u, \phi^p) \in V(\Omega).
    \end{align*}
    We will use this formulation for our numerical tests since it is equivalent to the $\dG(0)$ time discretization.
\end{remark}

\subsection{Weak space-time form of the adjoint problem}

The MORe DWR method \cite{FiRoWiChaFau2023} measures the adjoint sensitivity of the primal reduced-order solution with respect to some quantities of interest. Let a goal functional $J: \tilde{X}(\mathcal{T}_k, V(\Omega)) \rightarrow \mathbb{R}$ of the form
\begin{align}\label{eq_Ju_full}
    J(U) = \int_0^T J_1(U(t))\ \mathrm{d}t + J_2(U(T)),
\end{align}
be given, which represents some physical quantity of interest (QoI). The adjoint equations for porous media are derived, following \cite{Wi23_st}.
Since the problem is linear, we just need to switch trial and test functions in the bilinear form and replace the primal with the adjoint solution. The forcing function on the right-hand side is being substituted by the goal functional $J$, which is justified by the Lagrangian formalism, c.f. Section~\ref{sec_MOReDWR}. The space-time formulation of the adjoint problem thus reads:
Find $Z \in \tilde{X}\left(\mathcal{T}_k, V(\Omega)\right)$ such that
\begin{align*}
    A(\Phi)(Z)  = J^\prime_U(U)(\Phi) \qquad \forall \Phi \in \tilde{X}\left(\mathcal{T}_k, V(\Omega)\right).
\end{align*}
More concretely, the problem reads
\begin{align*}
    A(\Phi)(Z) &= \sum_{m = 1}^M \int_{I_m}  c(\partial_t \phi^p, z^p) + \alpha (\partial_t (\nabla_x \cdot \phi^u), z^p) + \frac{K}{\nu}(\nabla_x \phi^p, \nabla_x z^p) \ \mathrm{d}t \\
    &+ \sum_{m = 1}^M \int_{I_m} (\sigma(\phi^u), \nabla_x z^u) - \alpha(\phi^pI, \nabla_x z^u) + \alpha\langle \phi^pn, z^u\rangle_{\Gamma_N} \ \mathrm{d}t \\
    &+ \sum_{m = 1}^{M-1}\alpha([\nabla_x \cdot \phi^u]_{m}, z_{m}^{p,+}) + \alpha (\nabla_x \cdot \phi_{0}^{u,+},z_0^{p,+})
    + \sum_{m = 1}^{M-1}c([\phi^p]_{m}, z_{m}^{p,+}) + c (\phi_{0}^{p,+},z_0^{p,+}),
\end{align*}
where
\begin{align*}
    Z := (z^u, z^p), \quad \Phi := (\phi^u, \phi^p).
\end{align*}
The time derivatives are moved from the test function to the adjoint solution by integration by parts in time
e.g., on the time-continuous level we use integration by parts for 
\begin{align*}
    \int_0^T \alpha (\partial_t (\nabla_x \cdot \phi^u), z^p) = - \int_0^T \alpha (\nabla_x \cdot \phi^u, \partial_t z^p)\ \mathrm{d}t + \alpha (\nabla_x \cdot \phi^u(T), z^p(T)) - \alpha (\nabla_x \cdot \phi^u(0), z^p(0)).
\end{align*}
The adjoint problem runs backward in time and starts with a final time condition, which depends on the quantity of interest.
In particular, we have
\begin{align*}
    \alpha (\nabla_x \cdot \phi^u(T), z^p(T)) + c (\phi^p(T), z^p(T)) = J^\prime_{2, U}(U(T))(\Phi(T))
\end{align*}
i.e., $Z(T) = 0$ for $J_2 = 0$.
The integration by parts can be applied to each temporal element, which leads exemplarily to
\begin{align*}
    &\sum_{m = 1}^{M} \int_{I_m}  \alpha (\partial_t (\nabla_x \cdot \phi^u), z^p)  \ \mathrm{d}t +
     \sum_{m = 1}^{M-1}\alpha([\nabla_x \cdot \phi^u]_{m}, z_{m}^{p,+}) + \alpha (\nabla_x \cdot \phi_{0}^{u,+},z_0^{p,+}) \\
     &\qquad\quad = -\sum_{m = 1}^{M} \int_{I_m}  \alpha (\nabla_x \cdot \phi^u,\partial_t z^p)  \ \mathrm{d}t +
     \sum_{m = 1}^{M-1}\alpha([\nabla_x \cdot \phi^u]_{m}, z_{m}^{p,+}) + \alpha (\nabla_x \cdot \phi_{0}^{u,+},z_0^{p,+}) \\ 
     &\qquad\qquad\qquad\quad+ \sum_{m = 1}^{M}\alpha(\nabla_x \cdot \phi^{u,-}_{m}, z_{m}^{p,-}) - \alpha(\nabla_x \cdot \phi^{u,+}_{m-1}, z_{m-1}^{p,+}) \\
     &\qquad\quad = -\sum_{m = 1}^{M} \int_{I_m}  \alpha (\nabla_x \cdot \phi^u,\partial_t z^p)  \ \mathrm{d}t +
     \sum_{m = 1}^{M-1} \Big( \cancel{\alpha(\nabla_x \cdot \phi^{u,+}_{m}, z_{m}^{p,+})} - \alpha(\nabla_x \cdot \phi^{u,-}_{m}, z_{m}^{p,+})  \\ 
     &\qquad\qquad\qquad\qquad\quad+ \alpha(\nabla_x \cdot \phi^{u,-}_{m}, z_{m}^{p,-}) - \cancel{\alpha(\nabla_x \cdot \phi^{u,+}_{m}, z_{m}^{p,+})} \Big) \\
     &\qquad\qquad\qquad\quad  +\alpha (\nabla_x \cdot \phi_{M}^{u,-},z_{M}^{p,-}) - \bcancel{\alpha (\nabla_x \cdot \phi_{0}^{u,+},z_0^{p,+})} + \bcancel{\alpha (\nabla_x \cdot \phi_{0}^{u,+},z_0^{p,+})} \\
     &\qquad\quad = -\sum_{m = 1}^{M} \int_{I_m}  \alpha (\nabla_x \cdot \phi^u,\partial_t z^p)  \ \mathrm{d}t -
     \sum_{m = 1}^{M-1}\alpha(\nabla_x \cdot \phi^{u,-}_{m}, [z^{p}]_{m}) + \alpha (\nabla_x \cdot \phi_{M}^{u,-},z_{M}^{p,-}).
\end{align*}
Overall, we can rewrite the left-hand side of the adjoint problem as
\begin{equation}
\label{eq_adjoint_specific_equation_poroelasticity}
    \begin{aligned}
    A(\Phi)(Z) &= \sum_{m = 1}^{M} \int_{I_m}  -c(\phi^p, \partial_t z^p) - \alpha (\nabla_x \cdot \phi^u, \partial_t z^p) + \frac{K}{\nu}(\nabla_x \phi^p, \nabla_x z^p) \ \mathrm{d}t \\
    &+ \sum_{m = 1}^{M} \int_{I_m} (\sigma(\phi^u), \nabla_x z^u) - \alpha(\phi^pI, \nabla_x z^u) + \alpha\langle \phi^pn, z^u\rangle_{\Gamma_N} \ \mathrm{d}t \\
    &- \sum_{m = 1}^{M-1}\alpha(\nabla_x \cdot \phi^{u,-}_{m}, [z^{p}]_{m}) + \alpha (\nabla_x \cdot \phi_{M}^{u,-},z_{M}^{p,-})
    - \sum_{m = 1}^{M-1}c(\phi_{m}^{p,-}, [z^{p}]_{m}) + c (\phi_{M}^{p,-},z_{M}^{p,-}).
\end{aligned}
\end{equation}

\subsection{dG(0), i.e., backward Euler, time discretization of the adjoint problem}\label{sec:backward_euler_dual}
A $\dG(0)$ time discretization is again used for the adjoint problem.
Then on each temporal element $I_{m+1} = (t_{m}, t_{m+1})$, we have temporally constant functions $(z^u_m, z^p_m) \in V(\Omega)$ such that
\begin{align*}
    z^u\restrict{I_{m+1}} =: z^u_m, \qquad z^p\restrict{I_{m+1}} =: z^p_m. 
\end{align*}
Using piecewise-constant functions $f$ in time, the relations
\begin{align*}
    \partial_t f\restrict{I_m} = 0, \qquad -[f]_m = -(f_{m+1} - f_m) = f_m - f_{m+1}
\end{align*}
hold and give the adjoint time-stepping scheme: Find $Z_m := (z^u_m, z^p_m) \in V(\Omega)$ such that
\begin{align*}
    &c(z^{p}_m-z^{p}_{m+1}, \phi^p) + \alpha(z^{p}_m-z^{p}_{m+1}, \nabla_x \cdot \phi^u) + k\frac{K}{\nu}(\nabla_x z^{p}_m, \nabla_x \phi^p) \\
    &+ k(\nabla_x z^{u}_m, \sigma(\phi^u)) - k\alpha(\nabla_x z^{u}_m, \phi^pI) + k\alpha\langle z^{u}_m, \phi^p n\rangle_{\Gamma_N} = J^\prime_U(U)(\Phi)  \qquad \forall \Phi = (\phi^u, \phi^p) \in V(\Omega).
\end{align*}

\begin{remark}
    To be consistent with the backward Euler time discretization in Remark~\ref{remark:backward_euler_primal}, we again get rid of the timestep size  $k$ in
    displacement related terms and obtain
    \begin{align*}
        &c(z^{p}_m-z^{p}_{m+1}, \phi^p) + \alpha(z^{p}_m-z^{p}_{m+1}, \nabla_x \cdot \phi^u) + k\frac{K}{\nu}(\nabla_x z^{p}_m, \nabla_x \phi^p) \\
        &+ (\nabla_x z^{u}_m, \sigma(\phi^u)) - \alpha(\nabla_x z^{u}_m, \phi^pI) + \alpha\langle z^{u}_m, \phi^p n\rangle_{\Gamma_N} = J^\prime_U(U)(\Phi)  \qquad \forall \Phi = (\phi^u, \phi^p) \in V(\Omega).
    \end{align*}
\end{remark}

\section{MORe DWR: Model Order Reduction with Dual-Weighted Residual error estimates}
\label{sec_MOReDWR}

The MORe DWR approach aims at evaluating the quantities of interest using a reduced-order model, while guaranteeing the error due to this approximation. Thus, the difference between the reduced-order-model (ROM) solution $ U^{\text{ROM}} := U_{kh}^{\text{ROM}} \in X_k^{\dG(0)}(\mathcal{T}_k, V_h^{\text{ROM}})$ and the full-order-model (FOM) solution $U^{\text{FOM}} := U_{kh}^{\text{FOM}} \in X_k^{\dG(0)}(\mathcal{T}_k, V_h^{\text{FOM}})$ with $V_h^{\text{ROM}} \subset V_h^{\text{FOM}} =: V_h$ is controlled by employing a dual-weighted residual method \cite{BeRa96,becker_rannacher_2001,bangerth_rannacher_2003}. 
To this end, we obtain an optimization problem in which the model error between 
FOM and ROM, both measured in terms of some goal functional $J(\cdot)$, shall be minimized:
\begin{align}\label{eq:constrained_optimization_problem}
    J(U^{\text{FOM}}) - J(U^{\text{ROM}})
\end{align}
 subject to the constraint that the variational formulation of the poroelasticity problem (\ref{eq:space_time_form_poroelasticity}) is satisfied by $U^{\text{FOM}}$ i.e., $A(U^{\text{FOM}})(\Phi) = F(\Phi)$ for all test functions $\Phi \in X_k^{\dG(0)}(\mathcal{T}_k, V_h^{\text{FOM}})$.
For more information on space-time error control, we refer the reader to 
\cite{Schmich2009, ThiWi22_arxiv, RoThiKoeWi2023,Wi23_st}.
We focus in this work on the enrichment of the reduced basis depending on the temporal evolution of the quantities of interest. The reduced basis is refined in a goal-oriented way 
to accurately and efficiently compute the solution over the whole temporal domain.
In principle coarsening would also be possible, but
is not the objective in this work. For coarsening, we would need to follow 
the work of Meyer and Matthies \cite{Meyer2003}.

\subsection{Space-time dual-weighted residual method}
\label{sec:ST-DWR}

The two Lagrange functionals for the constrained optimization problem (\ref{eq:constrained_optimization_problem}) are defined as
\begin{align*}
      \mathcal{L}_\Box: X_k^{\dG(0)}(\mathcal{T}_k,V_h^\Box) \times X_k^{\dG(0)}(\mathcal{T}_k,V_h^\Box) & \rightarrow \mathbb{R}, \\
      (U^\Box , Z^\Box) & \mapsto J(U^\Box) - A(U^\Box)(Z^\Box) + F(Z^\Box)
\end{align*}
with $\Box \in \{\text{FOM}, \text{ROM}\}$.
The stationary points $(U^{\text{FOM}} , Z^{\text{FOM}} )$ and $(U^{\text{ROM}}, Z^{\text{ROM}})$ of the Lagrange functionals $\mathcal{L}_{\text{FOM}}$ and $\mathcal{L}_{\text{ROM}}$ need to satisfy the Karush-Kuhn-Tucker first-order optimality conditions. 
\subsubsection{Primal problem}
Firstly, the stationary points are solution to the following primal problems
\begin{align*}
    \mathcal{L}^\prime_{\Box,Z}(U^\Box , Z^\Box )(\delta Z^\Box) = -A(U^\Box)(\delta Z^\Box) + F(\delta Z^\Box) = 0 \qquad \forall \delta Z^\Box \in X_k^{\dG(0)},\quad \Box \in \{\text{FOM}, \text{ROM}\}.
\end{align*}
The corresponding solutions denoted $U^{\text{FOM}}$ and $U^{\text{ROM}}$ are the primal solutions. We observe that the primal solution can be obtained by solving the original problem, i.e., the quasi-static poroelasticity problem, forward in time as presented in Section~\ref{sec:backward_euler_primal}.

\subsubsection{Adjoint problem} \label{sec:adjoint_problem}
Secondly, the stationary points must also satisfy the following adjoint (or dual) problems
\begin{align*}
    \mathcal{L}^\prime_{\Box,U}(U^\Box , Z^\Box)(\delta U^\Box) = J^\prime_{U}(U^\Box)(\delta U^\Box)-A(\delta U^\Box)(Z^\Box) & = 0 \\ \forall \delta U^\Box &\in X_k^{\dG(0)}(\mathcal{T}_k,V_h^\Box), \quad \Box \in \{\text{FOM},\text{ROM}\},
\end{align*}
which give the adjoint solutions $Z^{\text{fine}}$ and $Z^{\text{coarse}}$. Note that here we already use that the poroelasticity equations are linear, which ensures that
\begin{align*}
    A^\prime_{U}(U)(\delta U, Z) = A(\delta U)(Z).
\end{align*}
Hence, the adjoint solution is obtained by solving
\begin{align}\label{eq:space_time_general_adjoint}
    A(\delta U)(Z) = J^\prime_{U}(U)(\delta U),
\end{align}
which is backward in time as discussed in Section~\ref{sec:backward_euler_dual}.

\begin{remark}\label{remark:adjoint_linear_problem}
   For linear goal functionals, the right-hand side of the adjoint problem (\ref{eq:space_time_general_adjoint}) reduces to
    \begin{align*}
        J^\prime_{U}(U)(\delta U) = J(\delta  U)
    \end{align*}
    and the adjoint problem does not depend on the primal solution anymore.
\end{remark}

\subsection{Error estimator}\label{sec:error_estimator}

The computable ROM error estimator is
\begin{align}\label{eq:computable_error_estimator}
    J(U^{\text{FOM}}) - J(U^{\text{ROM}}) \approx - A(U^{\text{ROM}})(Z^{\text{ROM}}) + F(Z^{\text{ROM}}) =: \eta,
\end{align}
following Theorem 4.2 from \cite{FiRoWiChaFau2023}. It is assembled on each temporal element $I_m$ separately to localize the error in time. More concretely, for poroelasticity the temporally-localized error estimator reads:
\begin{proposition}
Given the primal problem \eqref{eq:space_time_form_poroelasticity}, the  
goal functional \eqref{eq_Ju_full}, and the corresponding adjoint problem 
\eqref{eq:space_time_general_adjoint} with the adjoint bilinear form \eqref{eq_adjoint_specific_equation_poroelasticity},
the temporally-localized error estimator is
\begin{align*}
    \eta_m := \eta\restrict{I_m} =\, &-c(p_m-p_{m-1}, z^{p}_{m-1}) - \alpha(\nabla_x \cdot u_m-\nabla_x \cdot u_{m-1}, z^{p}_{m-1}) - k \cdot \frac{K}{\nu}(\nabla_x p_m,\nabla_x z^{p}_{m-1}) \\
	  &- (\sigma(u_m), \nabla_x z^{u}_{m-1}) + \alpha(p_mI, \nabla_x z^{u}_{m-1}) - \alpha\langle p_mn, z^{u}_{m-1}\rangle_{\Gamma_N} + \langle t, z^{u}_{m-1}\rangle_{\Gamma_N},
\end{align*}
for $m=1,\ldots,M$. Moreover, from \eqref{eq:computable_error_estimator}, we have
\[
J(U^{\text{FOM}}) - J(U^{\text{ROM}}) \approx \eta := \sum_{m=1}^M \eta_m. 
\]
\end{proposition}
The superscripts indicating reduced-order-model solutions for both the primal and dual solutions are omitted for a clearer presentation. The (reciprocal) effectivity index \cite{BaRhei78b}, 
which is the ratio between the estimated and the true errors, i.e., 
\begin{align}\label{eq:effectivity_index}
    \Ieff := \left|\frac{J(U^{\text{FOM}}) - J(U^{\text{ROM}})}{\eta}\right|
\end{align}
is used to measure the quality of our error estimator. We desire $\Ieff \approx 1$  since then the error estimator can reliably predict the reduced-order-modeling error. We refer to \cite{EndtLaWi20} for two-sided proofs of discretization errors measured in goal functionals using saturation assumptions. We also observe that the effectivity index $\Ieff$ is close to one in the numerical tests in Section~\ref{sec:numerical_tests}.
Finally, the quality of the adaptive refinement is measured in terms of the indicator 
index \cite{RiWi15_dwr}:
\begin{align}\label{eq:indicator_index}
    \Iind := \frac{|J(U^{\text{FOM}}) - J(U^{\text{ROM}})|}{\sum_{m=1}^M |\eta_m|}.
\end{align}

\subsection{Error estimator-based ROM updates}

In this section, we give a compact summary of 
our novel approach of a 
goal-oriented incremental reduced-order model from \cite{FiRoWiChaFau2023} and discuss its extension 
to porous media settings.
In the MORe DWR method, we marry a reduced-order model with a dual-weighted-residual-based error estimator and an incremental version of the POD algorithm. The goal is then to use the error estimator to identify when the solution behavior is not captured accurately by the reduced basis, such that it is incrementally enriched on-the-fly with new full-order-model snapshots. The approach is particularly cheap because the FOM snapshots are computed only for one temporal element.

In more detail, we apply our findings on error control of Section~\ref{sec:ST-DWR} to a backward Euler reduced-order model of poroelasticity. Further, an incremental basis generation is mandatory for the method to reduce computational operations and thus to be fast, which we realize by means of an incremental SVD. The incremental SVD is presented in Section~\ref{sec:additive_svd}. In this context, we also introduce the incremental POD as a trimmed version of the incremental SVD.
Subsequently, the overall MORe DWR framework is depicted in Section~\ref{sec:incremental_ROM}, where all the ingredients are assembled and the final algorithm is presented.

In summary, our novel approach avoids a computationally heavy offline phase and directly solves the reduced model. Our approach can be thought of as a replacement for traditional full-order-model simulations, since we highly reduce the computational cost, while controlling the error between the full-order (FOM) and the reduced-order model (ROM).

\subsubsection{Incremental Proper Orthogonal Decomposition}
\label{sec:additive_svd}

The Proper Orthogonal Decomposition (POD) \cite{lassila_model_2014,benner2020model,kunisch2002galerkin,ravindran2000reduced,sirovich1987turbulence,willcox2002balanced,christensen1999evaluation,gunzburger2007reduced,caiazzo2014numerical,baiges2013explicit, Ekre2020} can be employed to project the differential equation onto a lower-dimensional solution manifold. In traditional applications of POD, one first collects full-order solution snapshots
\begin{align*}
    U(t_0), \, U(t_1),\,\, \dots,\,\, U(t_M),
\end{align*}
and finds the POD basis by computing the singular value decomposition
\begin{align*}
    Y := \begin{bmatrix}
        U(t_0) & U(t_1) & \dots & U(t_M)
    \end{bmatrix} = \Psi \Sigma \Phi^T
\end{align*}
with orthogonal matrices $\Psi \in \mathbb{R}^{n \times d}, \Phi \in \mathbb{R}^{(M+1) \times d}$ and $\Sigma = \operatorname{diag}(\sigma_1, \dots, \sigma_d) \in \mathbb{R}^{d\times d}$. Here, $n$ is the number of spatial degrees of freedom, $M+1$ is the number of 
FOM snapshots, and $d \leq \min(n,M+1).$
The SVD is then being truncated to the size $N \ll d$ 
using a retained energy criterion, cf. \cite{grassle2018pod, gubisch2017proper, lassila_model_2014} e.g.,
\begin{align}
    \varepsilon(N) = \frac{\sum_\mathrm{i=1}^{N} \sigma_{i}^2}{\sum_{i=1}^{d} \sigma_i^2} = \frac{\sum_{i=1}^{N} \sigma_i^2}{\sum_{i=1}^\mathrm{M+1} ||{U}_i||_{\mathbb{R}^n}^2} < 99.9 \%.
    \label{eq:energy_content}
\end{align}
For 
coupled problems with several PDEs, like poroelasticity, we could again perform the SVD of the snapshot matrix $Y$ consisting of the entire solution vectors $U(t_i) = \begin{pmatrix}
    u(t_i) \\ p(t_i)
\end{pmatrix}$, but we choose to create separate SVDs for displacement $u(t_i)$ and pressure $p(t_i)$ such that POD bases with various sizes can be used for the different solution components.

\begin{figure}[H]
    \begin{center}
    \includegraphics[width=\textwidth]{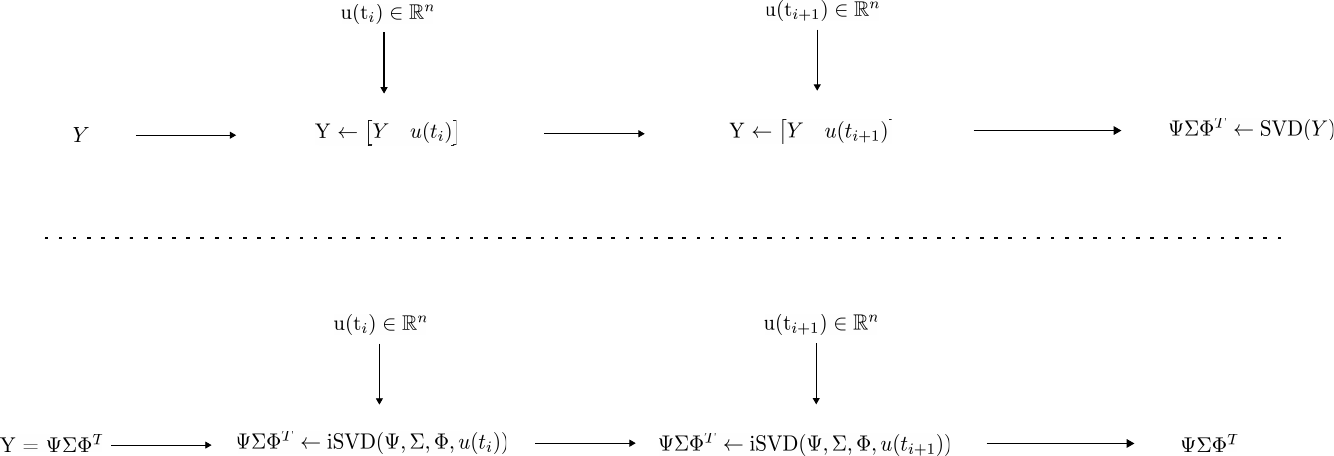}
    \caption{Methodology of SVD (top) and incremental SVD (bottom).}
    \label{fig:comparison_svd_isvd}
    \end{center}
\end{figure}

Instead of building repeatedly a reduced basis from scratch, we suggest updating an already existing truncated SVD (tSVD) or solely its left-singular (POD) vectors according to modifications of the snapshot matrix without recomputing the whole tSVD or requiring access to the snapshot matrix \cite{kuhl2023incremental, FiRoWiChaFau2023}.
The POD becomes incremental by appending additional snapshots to the initial snapshot matrix. This difference between the classical SVD and incremental SVD (iSVD) is illustrated in Figure~\ref{fig:comparison_svd_isvd}.  In this context, we rely on the general approach of an additive rank-$b$ modification of the SVD, mainly developed by \cite{brand2002incremental, brand2006fast} and applied to the model-order reduction of fluid flows in \cite{kuhl2023incremental}. 

Let ${Y} \in \mathbb{R}^{n \times \mathrm{\Tilde{m}}}$ be a given snapshot matrix that includes $\Tilde{m}>0$ snapshots. Usually, $\Tilde{m}$ is equal or connected to the number of already computed time steps. Further, we have the rank-$N$ tSVD $\Psi \Sigma \Phi^\mathrm{T}$ of the matrix $Y$. Additionally, let $b\in\mathbb{N}$ newly computed snapshots $\{U_1, \dots, U_b\}$ be stored in the bunch matrix
\begin{align} \label{eq:bunch_matrix}
    B = \begin{bmatrix}
            u_1 & \dots & u_b
        \end{bmatrix} \in \R^{n \times b}.
\end{align}
We aim to compute the tSVD that is updated by the information contained in the bunch matrix $B$ according to
\begin{align*}
    \Tilde{\Psi} \Tilde{\Sigma} \Tilde{\Phi}^{T} =
    \Tilde{Y} =
    \begin{bmatrix}
        Y & B
    \end{bmatrix}
\end{align*}
without explicitly recomputing $Y$ or $\Tilde{Y}$ due to performance and memory reasons which was the original motivation of Brand's work on the incremental SVD, cf. \cite{brand2002incremental,brand2006fast}.
We skip the technical derivations of the incremental SVD in this work and instead summarize the incremental POD routine in Algorithm~\ref{algo:iPOD}. A prototypical implementation of the incremental POD in Python can be found at \url{https://github.com/Hendrik240298/Incremental_POD} and for a mathematical deep dive into the theory behind the incremental POD, we refer the reader to \cite{kuhl2023incremental}[Section 2.2] and \cite{FiRoWiChaFau2023}[Section 4.2.1].
\begin{algorithm}[H]
    \caption{Incremental POD update} \label{algo:iPOD}
    \hspace*{\algorithmicindent} \textbf{Input:} Reduced basis matrix $\Psi_N \in \R^{n \times N}$,
    singular value vector ${\Sigma = [\sigma_1, \dots, \sigma_N] \in \R^N}$,
    bunch matrix $B \in \R^{n \times b}$,
    and energy threshold $\varepsilon \in [0,1]$.\\
    \hspace*{\algorithmicindent} \textbf{Output:} Reduced basis matrix $\Tilde{\Psi}_N \in \R^{n \times \Tilde{N}}$,
    singular value vector $\Tilde{\Sigma} = [\Tilde{\sigma}_1, \dots, \Tilde{\sigma}_{\Tilde{N}}]$
    \begin{algorithmic}[1]
        \State $H = \Psi_N^T B$
        \State $P = B - \Psi_N H$
        \State $[Q_P, \, R_P] = \text{QR}(P)$
        \State $Q = [\Psi_N \; Q_P] $
        \State ${F} = \begin{bmatrix}
                {\text{diag}(\Sigma)} & {H}     \\
                {0}                   & {R}_{P}
            \end{bmatrix}$
        \If{Q not orthogonal}
        \State $[Q,\, R] = \text{QR}(Q)$
        \State $ F = RF$
        \EndIf
        \State $[\Psi^\prime, \Sigma^\prime] = \text{SVD}(F) $
        \State $\Tilde{N} = \min\left\lbrace N \in \mathbb{N}\ \middle|\ \varepsilon(N) \geq \varepsilon,\ \ 1 \leq N \leq d \right\rbrace$
        \State $\Tilde{\Sigma} = \text{diag}(\Sigma^\prime)(1:\Tilde{N})$
        \State $\Tilde{\Psi}_N = Q \Psi^\prime(:,1:\Tilde{N})$
    \end{algorithmic}
\end{algorithm}
\subsection{Goal-oriented error-controlled incremental ROM}
\label{sec:incremental_ROM}
This section outlines how goal-oriented error control (Section~\ref{sec:ST-DWR})
drives the incremental reduced-order model (Section~\ref{sec:additive_svd}).

\subsubsection{The MORe DWR algorithm}
\begin{figure}[H]
    \centering
    \resizebox{0.7\textwidth}{!}{
        \begin{tikzpicture}[node distance=2.2cm]
          \node[rounded corners, rectangle, draw, align=center] (pROM) {{\color{blue} \textbf{Primal ROM}} for $I_1, I_2, \dots, I_M$};
          \node[rounded corners, rectangle, draw, below left of=pROM,xshift=-2cm, align=center] (dROM) {{\color{purple} \textbf{Dual ROM}} for $I_M, I_{M-1} \dots, I_1$};
          \node[rounded corners, rectangle, draw, below of=dROM, align=center] (esti) {{\color{green!50!black} \textbf{Error estimates}} for $\eta_{1}, \eta_{2}, \dots, \eta_{M}$ \\ \textbf{Localize} $m_{\max} := \underset{1 \leq m \leq M}{\argmax}\,\, \eta_{m}^{\text{rel}}$};
          \node[diamond, aspect=1.5, draw, below right of=esti, xshift=2cm, yshift=-0.5cm, align=center] (desc) {\textbf{$\eta^{\text{rel}} < \mathrm{tol}$}};
          \node[rounded corners, rectangle, draw, above right of=desc, xshift=2cm, yshift=0.25cm, align=center] (pFOM) {\textbf{Primal FOM on $I_{m_{max}}$} 
          \\ + \\ \textbf{Enrich primal RB}};
          \node[rounded corners, rectangle, draw, above of=pFOM, yshift=0.25cm, align=center] (dFOM) {\textbf{Dual FOM on $I_{m_{\max}}$}\\+\\ \textbf{Enrich dual RB}};
          \node[below right of=desc,yshift=0.5cm] (end) {Done};
          
          \draw[->, rounded corners=20pt] (pROM.west) -| (dROM.north);
          \draw[->] (dROM) -- (esti);
          \draw[->, rounded corners=20pt] (esti) |- (desc.west);
          \draw[->, rounded corners=20pt] (desc) -| node[pos=0.05, above, align=center] {No} (pFOM);
          \draw[->] (pFOM) -- (dFOM);
          \draw[->, rounded corners=20pt]  (dFOM) |- (pROM);
          \draw[->, rounded corners]  (desc) |- node[pos=0.85, above, align=center] {Yes} (end);
      \end{tikzpicture}
    }
    \caption{Schematic representation of MORe DWR algorithm.}
    \label{fig:moredwr_tikz_illustration}
\end{figure}
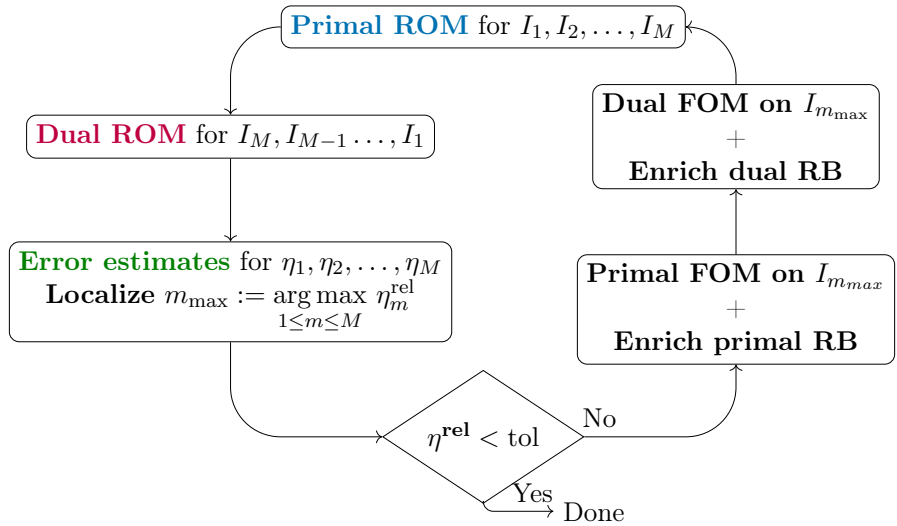

The overall procedure is illustrated in Figure~\ref{fig:moredwr_tikz_illustration}.

The aim is to solve the reduced-order model and adaptively enrich the reduced basis by means of the iPOD with full-order solutions until a given estimated error tolerance is reached for the 
chosen goal functional. The approach is designed to work without any prior knowledge or exploration of the solution manifold, while also attempting to minimize the number of full-order solves. The error associated with the reduced solution is estimated on each temporal element. 
The new full-order solution used for the basis enrichment is computed on the temporal element with the largest error. So, two full-order solves are conducted for each enrichment. We observe that this procedure is perfectly compatible with the adaptive basis selection based on DWR estimates presented by Meyer and Matthies in \cite{Meyer2003} to reduce the dimension of the reduced space. Thus, if incorporated it would be possible to either enrich or downsize the reduced basis according to the problem statement.
The resulting approach is detailed in Algorithm~\ref{algo:iROM}.
The reduced primal and dual solutions are used to provide a fast error estimation.
Thus, the incremental ROM has both a primal and a dual reduced basis as input in the form of the reduced basis matrix.

\subsubsection{Initialization}
If no prior information is available, the bases are initialized by computing the primal and dual solution snapshots on the first temporal element. If prior reduced bases are available e.g., from a previous simulation on a different parameter configuration or some (very cheap) coarse-grid solution in the context of a multigrid idea, prior bases are 
used as an initial guess and are then altered by the goal-oriented adaptation.
Thus, the MORe DWR approach is totally appropriate for reduced-order modeling of parameterized problems. 
The method can be an efficient substitute for the full-order model, as it builds a reduced basis tailored to a quantity of interest with a minimum of FOM solves.

Nevertheless, in the scenario without prior knowledge, we observed an issue when exclusively the first dual FOM solution ($t = T$) is utilized for the basis generation. We identified a significant discrepancy between the error estimates and the real errors, which is demonstrated and further described in Figure~\ref{fig:adaptivity_footing_without_dual}. 
To mitigate this issue, further assessments utilizing the final dual snapshots for the basis enrichment at the first temporal elements ($t = 0$) were conducted, which consistently improved the accuracy of the error estimates. Nonetheless, a practical constraint is posed, as these snapshots are not readily accessible without completing the whole dual FOM solution process.
Thus, a computationally efficient surrogate to bypass the need for the full-order solution process is required. For this, the dual FOM is substituted by the dual ROM, which is already computed during the enrichment process and only the last dual solution steps are conducted by the dual FOM. 
Given that the inconsistency in the estimate predominantly emerges during the initial enrichment iterations, we have optimized the process by confining the additional dual ROM enrichment solely to these early iterations for efficiency.

\subsubsection{Adaptive algorithm}
At each MORe DWR iteration, the primal and dual reduced bases are enriched for the temporal element $m_{max}$ corresponding to the highest relative error comparing the local relative error
\begin{equation}
    \eta_m^{rel} = \frac{\eta_m}{J(U^{\text{ROM}}) + \sum_{l = 1}^M \eta_l}
\end{equation}
associated with any temporal element $m$. Because the full-order solutions are not available, the error is normalized according to
\begin{equation}
        J(U^{\text{FOM}}) \approx J(U^{\text{ROM}}) + \eta = J(U^{\text{ROM}}) + \sum_{l = 1}^M \eta_l.
\end{equation}
The adaptive MORe DWR algorithm is stopped according to the global relative error estimator
\begin{equation}
    \eta^{\text{rel}} = \frac{\eta}{J(U^{\text{FOM}})} \approx \frac{\sum_{m=1}^M \eta_m}{J(U^{\text{ROM}}) + \sum_{m = 1}^M \eta_m}  < TOL^{rel}
\end{equation}
with $TOL^{rel}$ a user-chosen threshold. 
This ensures that relative error over the entire temporal domain is below a given tolerance, 
whereas in \cite{FiRoWiChaFau2023} we ensured that the local relative error is below a given tolerance. However, enforcing the accuracy locally is a more greedy approach, which leads to more full-order-model solves.
\begin{algorithm}[H]
    \caption{MORe DWR algorithm}
    \label{algo:iROM}
    \hspace*{\algorithmicindent} \textbf{Input:} Initial condition $U_0:=U(t_0)$, primal and dual reduced basis matrices $(\Psi^{p_u}_{N_p^u}, \Psi^{p_p}_{N_p^p})$ and $(\Psi^{d_u}_{N_d^u},\Psi^{d_p}_{N_d^p})$, energy iPOD threshold $\varepsilon \in [0,1]$, and error tolerance $\text{tol}>0$.\\
    \hspace*{\algorithmicindent} \textbf{Output:} Primal and dual reduced basis matrices  $(\Psi^{p_u}_{N_p^u}, \Psi^{p_p}_{N_p^p})$ and $(\Psi^{d_u}_{N_d^u},\Psi^{d_p}_{N_d^p})$ and reduced primal solutions $U_{m}$ for all $1\leq m \leq M$.
    \begin{algorithmic}[1]
        \While{$\eta^{\text{rel}} > tol$}
        \For{$m=1, 2, \dots, M$} \Comment{Primal ROM}
        \State Solve reduced primal system on $I_m$ for $U_m^{\text{ROM}}$
        \EndFor
        \For{$m=M, M-1, \dots, 1$} \Comment{Dual ROM}
        \State Solve reduced dual system on $I_m$ for $Z_{m-1}^{\text{ROM}}$
        \EndFor
        \State Compute error estimate: $\eta^{\text{rel}}$
        \If{$\eta^{\text{rel}} > tol$}
        \State ind temporal element with maximum error: ${m_{max}} =  \argmax\limits_{1 \leq m \leq L} \left| \frac{\eta_m}{J(U^{\text{ROM}}) + \sum_{l = 1}^M \eta_l} \right|$
        \vspace{0.2em}
        \State Solve primal full-order system on $I_{m_{max}}$ for $U_{m_{max}}^{\text{FOM}} = (u_{m_{max}}^{\text{FOM}}, p_{m_{max}}^{\text{FOM}})$
        \vspace{0.2em}
        \State Update primal reduced basis: 
        \begin{align*}
            \Psi^{p_u}_{N^u_p} &= \text{iPOD}(\Psi^{p_u}_{N^u_p}, \Sigma_{N^u_p}, [u_{m_{max}}^{\text{FOM}}], \varepsilon) \\
            \Psi^{p_p}_{N^p_p} &= \text{iPOD}(\Psi^{p_p}_{N^p_p}, \Sigma_{N^p_p}, [p_{m_{max}}^{\text{FOM}}], \varepsilon) 
        \end{align*} 
        \vspace{-0.5cm}
        \State Solve dual full-order system on $I_{m_{max}}$ for $Z_{m_{max}-1}^{\text{ROM}} = (z_{m_{max}-1}^{u,\text{FOM}}, z_{m_{max}-1}^{p,\text{FOM}})$
        \vspace{0.2em}
        \State Update dual reduced basis: 
        \begin{align*}
            \Psi^{d_u}_{N^u_d} &= \text{iPOD}(\Psi^{d_u}_{N^u_d}, \Sigma_{N^u_d}, [z_{m_{max}-1}^{u,\text{FOM}}], \varepsilon) \\
            \Psi^{d_p}_{N^p_d} &= \text{iPOD}(\Psi^{d_p}_{N^p_d}, \Sigma_{N^p_d}, [z_{m_{max}-1}^{p,\text{FOM}}], \varepsilon) 
        \end{align*} 
        \vspace{-0.5cm}
        \State Update reduced system components and error estimator
        \EndIf
        \EndWhile
    \end{algorithmic}
\end{algorithm}

\begin{remark}
    The additional enrichment of the dual space to provide an accurate estimation of the error is omitted in Algorithm~\ref{algo:iROM} and in Figure~\ref{fig:moredwr_tikz_illustration}  for the sake of clarity.  This dual enrichment takes place after line 12 in Algorithm~\ref{algo:iROM} and after enriching the primal and dual reduced bases in Figure~\ref{fig:moredwr_tikz_illustration}. The supplementary enrichment is based on the FOM solution of the dual problem for a few time steps, corresponding to the few last time steps of the dual problem i.e., the few first time steps of the primal problem. As the dual problem is backward in time, it starts from the ROM solution at $t = t_l$ with $l \in \mathbb{N}$ and the solution is computed for $t_{l-1}$, $t_{l-2}$, \dots, $t_0 $. The dual reduced basis is then enriched with these additional snapshots. This procedure is applied only in the few first MORe DWR iterations prior to the update of the reduced system components and the error estimator.
\end{remark}

\section{Numerical tests}
\label{sec:numerical_tests}

The MORe DWR framework is numerically substantiated on two different problem configurations. 
The Mandel problem 
\cite{Mandel1953,Cr63,Cheng88,AbChCuDeRo96,DuiMiWi22} (see also \cite{Gai2004, Guzman2012, Girault2011, Liu04, Wick2020PFF} and more recently also for nonlinear poroelasticity \cite{DuiMiWi22}) 
is first considered as a benchmark for two-dimensional poroelasticity and a linear goal functional. 
In Mandel's problem one observes the so-called Mandel-Cryer effect \cite{Cr63,DuiMi21} of 
a non-monotonic pressure evolution: first increasing pressure, followed by decreasing pressure. 
The second numerical test is a three-dimensional footing problem inspired by \cite{Gaspar2008}. 
Our computations have been performed on an AMD EPYC 7H12 64-Core Processor. 
The FEM codes have been written in FEniCS \cite{fenics2015} and
the reduced-order modeling has been performed with NumPy \cite{NumPy2020} and SciPy \cite{SciPy2020}. 

\subsection{Mandel problem in 2D}\label{sec:mandel}
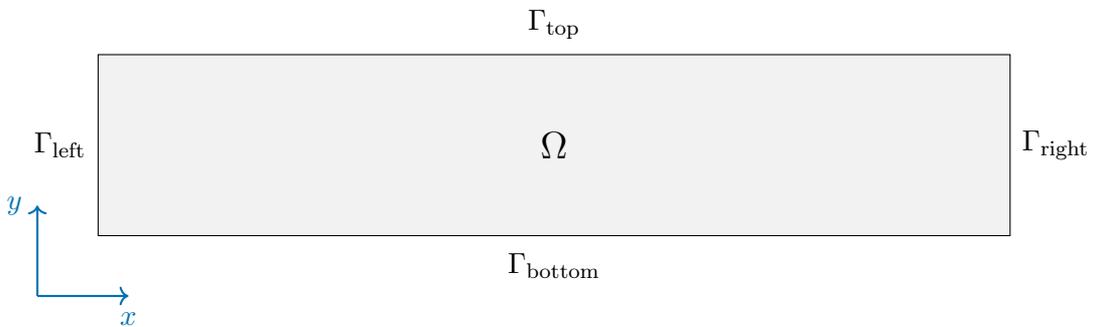
\begin{figure}[H]
    \begin{center}
    \begin{tikzpicture}[scale = 4, draw=black]
        \draw[draw=black, fill=black!5!white] (0,0) rectangle (3,0.6);
        \node (omega) at (1.5,0.3) {\Large{$\Omega$}};
        \node (left) at (-0.125,0.3) {$\Gamma_{\text{left}}$};
        \node (right) at (3.15,0.3) {$\Gamma_{\text{right}}$};
        \node (top) at (1.5,0.7) {$\Gamma_{\text{top}}$};
        \node (bottom) at (1.5,-0.1) {$\Gamma_{\text{bottom}}$};

        \draw[draw=blue, ->, thick]
        (-0.2,-0.2) -- (0.1,-0.2);
        \node(x) at (0.1,-0.275) {\color{blue}$x$};
        \draw[draw=blue, ->, thick]
        (-0.2,-0.2) -- (-0.2,0.1);
        \node(y) at (-0.275,0.1) {\color{blue}$y$};
    \end{tikzpicture}
    \caption{Domain for Mandel's problem.}
    \label{fig:mandel_domain}
    \end{center}
\end{figure}

Let $\Omega := (\SI{0}{\meter}, \SI{100}{\meter}) \times (\SI{0}{\meter}, \SI{20}{\meter})$ and $I := (\SI{0}{\second}, \SI{5000000}{\second}) \approx  (\SI{0},\SI{58}{\day})$ with boundaries as shown in Figure~\ref{fig:mandel_domain}.

\noindent The initial and boundary conditions are given by
\begin{align*}
    p(\cdot, 0) &= p^0(\cdot) = 0 &&\text{in } \Omega \times \{0\}, \\
    u(\cdot, 0) &= u^0(\cdot) = 0 &&\text{in } \Omega \times \{0\}, \\
    \frac{K}{\nu}\nabla_x p\cdot n &= 0 &&\text{on } \partial \Omega \setminus \Gamma_{\text{right}} \times I, \quad\tag*{(No flow condition, homogeneous Neumann)}\\
    \sigma(u) \cdot n &= -\bar{t}e_y &&\text{on } \Gamma_{\text{top}} \times I, \qquad\quad\text{(Traction condition, inhomogeneous Neumann)}\\
    p &= 0  &&\text{on } \Gamma_{\text{right}} \times I, \quad\tag*{(Constant zero pressure, homogeneous Dirichlet)}\\
    \sigma(u) \cdot n &= 0 &&\text{on } \Gamma_{\text{right}} \times I,\quad\tag*{(Traction free, homogeneous Neumann)}\\
    u_y &= 0 \quad\text{and}\quad \partial_y u_x = 0  &&\text{on } \Gamma_{\text{bottom}} \times I,\quad\text{(Confined conditions, mixed Dirichlet/Neumann)} \\
    u_x &= 0  \quad\text{and}\quad \partial_x u_y = 0 &&\text{on } \Gamma_{\text{left}} \times I, \quad\tag*{(Confined conditions, mixed Dirichlet/Neumann)}
\end{align*}
The parameters for Mandel's problem are summarized in Table~\ref{tab:params_mandel}.
\begin{table}[H]
    \centering
    \begin{tabular}{ |p{2.5cm}||p{2.5cm}|}
         \hline
         Parameter & Value \\
         \hline
            M & \SI{1.75e7}{\pascal} \\
            c & 1/M \\
            $\alpha$ & \SI{1}{\pascal\metre} \\
            $\nu$ & \SI{1e-3}{\metre\squared\per\second} \\
            K & \SI{1e-13}{\metre\squared} \\
            $\rho$ & \SI{1}{\kilogram\per\metre\cubed} \\
            $\bar{t}$ & \SI{1e7}{\pascal\metre} \\
            $\mu$ & \SI{1e8}{} \\
            $\lambda$ & $\frac{2}{3} \times 10^8$ \\
         \hline
    \end{tabular}
    \caption{Parameters in Mandel's problem.} 
    \label{tab:params_mandel}
\end{table}
\noindent In space, we use Taylor-Hood elements, i.e., quadratic finite elements for the displacements $u$ and linear finite elements for the pressure $p$. 
For this numerical test, a fixed spatial mesh with 80 spatial cells in the $x$-direction and 16 spatial cells in the $y$-direction is designed, which leads to an isotropic mesh with $10,626$ DoFs for displacement and $1,377$ DoFs for pressure.
For the temporal discretization, the end time is $T = \SI{5000000}{\second} \,  \approx \SI{58}{\day}$ with $5,000$ temporal elements i.e., the chosen time step size $k = \SI{1000}{\second}$. For the quantity of interest, we choose the time-integrated pressure acting on the bottom boundary i.e., 
\begin{align*}
    J(U) := \int_I \int_{\Gamma_{\text{bottom}}} p\ \mathrm{d}x\ \mathrm{d}t.
\end{align*}
Further, for the MORe DWR enrichment, we choose $1-10^{-7}$ and $1-10^{-11}$ for the primal displacement and pressure energy tolerances \eqref{eq:energy_content}. Accordingly, the dual energy tolerances are both set to $1-10^{-9}$. As argued in Section~\ref{sec:incremental_ROM}, it is crucial to have a dual space large enough for estimating accurately the error. Thus, the reduced dual basis is enriched during the first MORe DWR loops from the snapshots of the full-order dual solution on the first time steps $[t_0,t_{\ell -1}]$ i.e., for the last time points of the dual problem which is formulated backward in time. Here, it is enriched for the first five MORe DWR iterations based on the snapshots corresponding with $[t_0,t_4]$ in the primal temporal discretization.

\begin{figure}[H]
    \centering
    \begin{subfigure}[b]{0.485\textwidth}
        \centering
        \includegraphics[width=\textwidth]{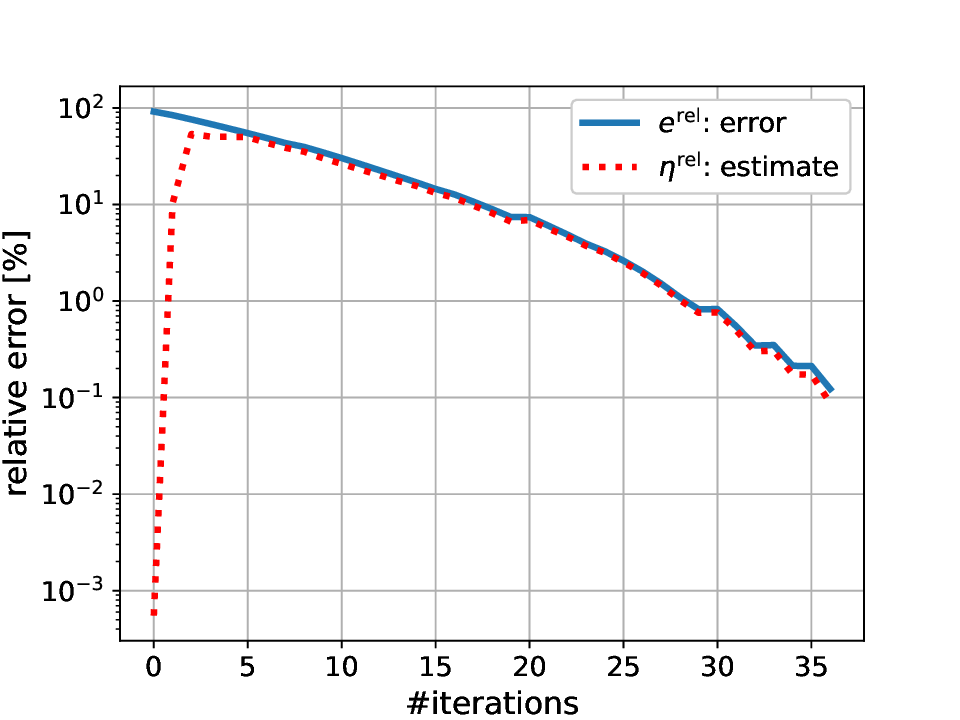}
        \caption{Evolution of the true and estimated relative errors throughout the MORe DWR iterations.}
        \label{fig:adaptivity_mandel:iteration_error}
    \end{subfigure}
    \hfill
    \begin{subfigure}[b]{0.485\textwidth}
        \centering
        \includegraphics[width=\textwidth]{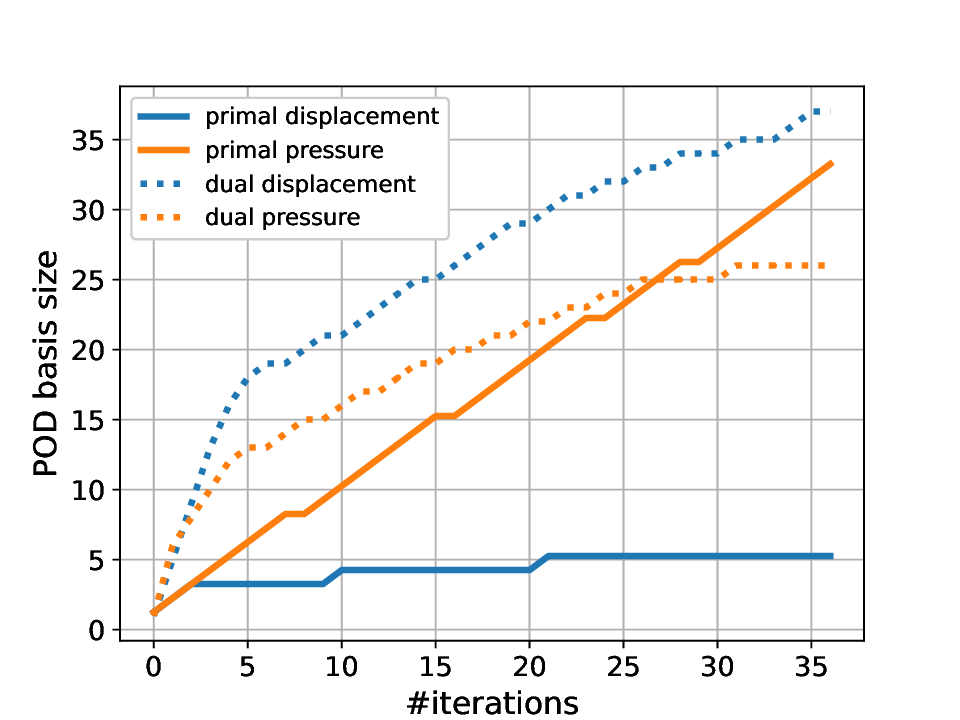}
        \caption{Evolution of the POD basis sizes throughout the MORe DWR iterations. \phantom{placeholder}}
        \label{fig:adaptivity_mandel:basis_size}
    \end{subfigure}
    \caption{Adaptivity in the MORe DWR method of Mandel's problem.}
    \label{fig:adaptivity_mandel}
\end{figure}

First, the adaptive enrichments of the method are shown in Figure~\ref{fig:adaptivity_mandel}. The true relative error between the reference full-order solution and the reduced-order solution given by
\begin{align*}
    e^{rel} := \frac{|J(U^{\text{FOM}})-J(U^{\text{ROM}})|}{|J(U^{\text{FOM}})|}.
\end{align*} 
Here, we observe a clearl decrease with the iterations of the MORe DWR algorithm; see Figure~\ref{fig:adaptivity_mandel:iteration_error}. In contrast 
to the original MORe DWR paper \cite{FiRoWiChaFau2023},
the relative error is now much closer to or even 
slightly exceeds the error tolerance. This 
can be explained by the fact that we do not 
enforce the error tolerance on each temporal element,
but only on the whole integrated time domain.
The estimated relative error is very close to the true error after iteration $2$. Thereafter, only slight deviations can be noticed locally. On the contrary, a large discrepancy can be observed between the true and estimated errors for the first iterations, where the estimate yields significant underestimations due to the inaccuracy of the dual reduced space. The target tolerance of $0.1 \%$ is reached after 36 iterations. For larger error tolerances, the algorithm would terminate sooner e.g., for $1 \%$ relative error, we need 29 iterations, for $5 \%$ error we need 22 iterations and for $20 \%$ error only 13 iterations are required (see Table~\ref{tab:comparison:mandel}).
The enrichment of the POD basis is illustrated in Figure~\ref{fig:adaptivity_mandel:basis_size}. The primal reduced spaces increases progressively and regularly with the iterations. We observe that the primal pressure solution requires significantly more modes than the primal displacement. 
This implies that the choice of separating the bases for displacement and pressure is beneficial. Interestingly, the dual displacement solution (37 modes) requires more POD modes than the primal displacement solution (5 modes), whereas the dual pressure solution (26 modes) requires less POD modes than the primal pressure solution (33 modes). We also see that the sizes of the dual bases increase faster than in the primal case during the first 5 iterations as we enforce an additional dual basis enrichment. But, its slope decreases after these initial enrichments.

\begin{figure}[H]
    \centering
    \begin{subfigure}[b]{0.485\textwidth}
        \centering
        \includegraphics[width=\textwidth]{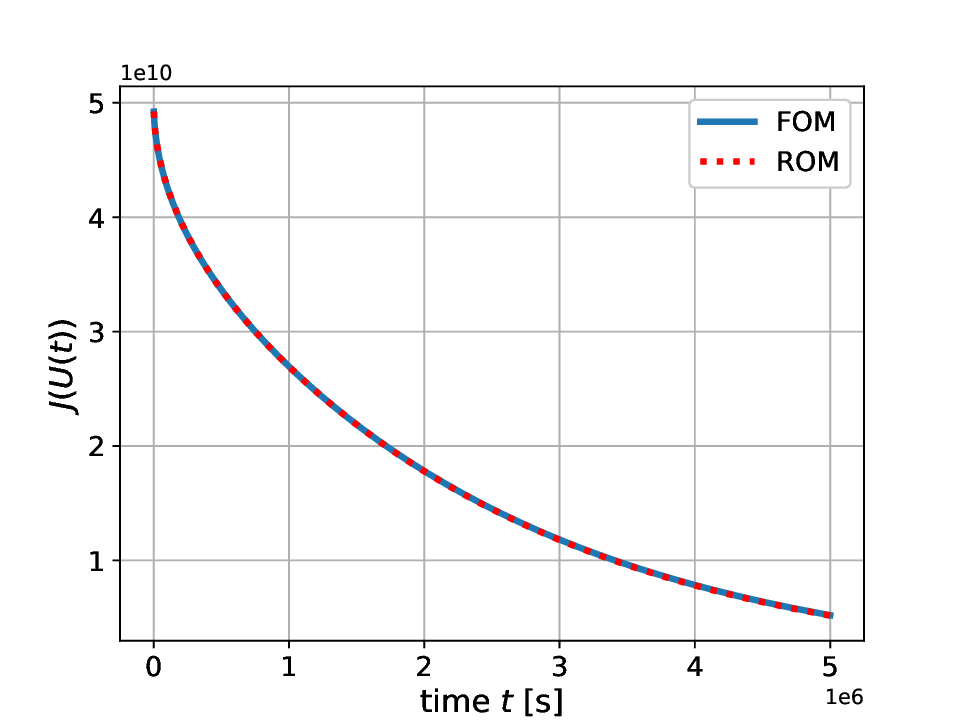}
        \caption{$TOL^{rel} = 0.1\%$}
    \end{subfigure}
    \hfill
    \begin{subfigure}[b]{0.485\textwidth}
        \centering
        \includegraphics[width=\textwidth]{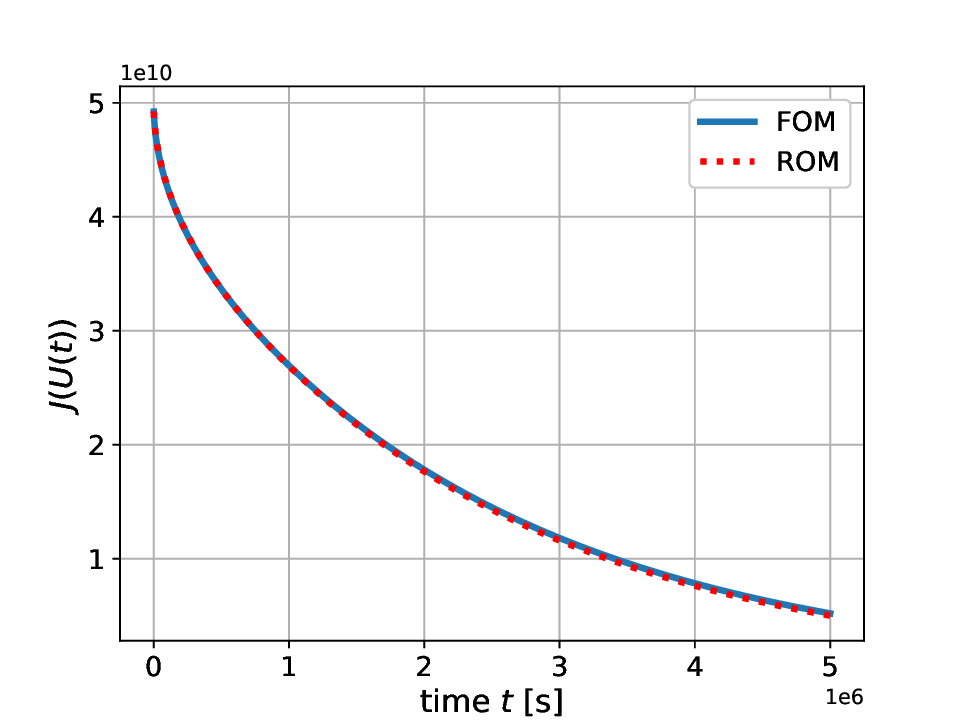}
        \caption{$TOL^{rel} = 1.0\%$}
    \end{subfigure}
        \begin{subfigure}[b]{0.485\textwidth}
        \centering
        \includegraphics[width=\textwidth]{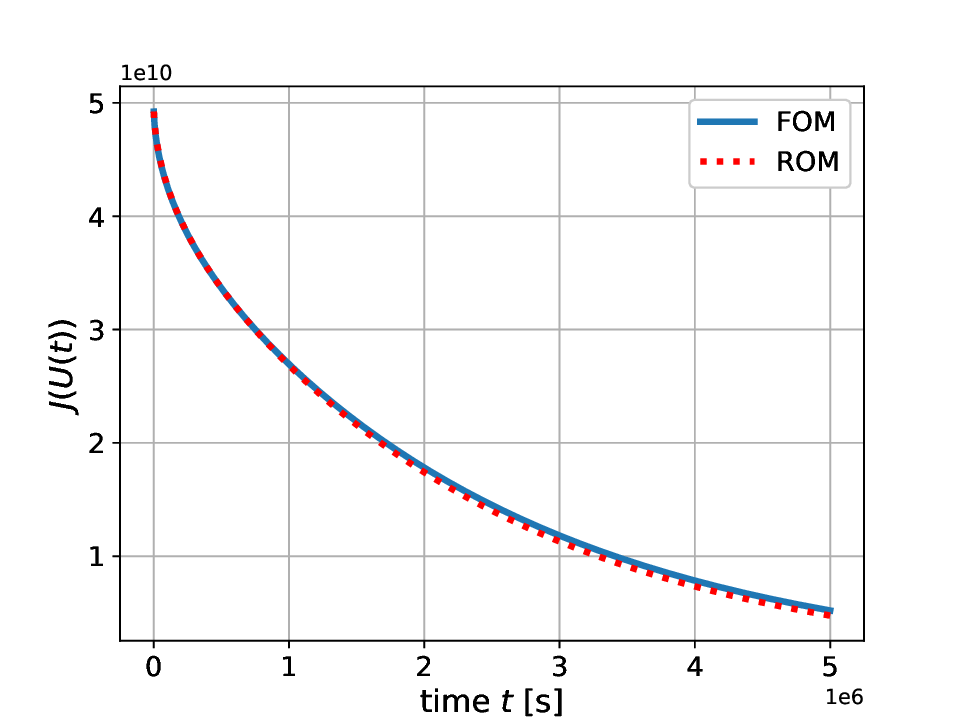}
        \caption{$TOL^{rel} = 2.0\%$}
    \end{subfigure}
    \begin{subfigure}[b]{0.485\textwidth}
        \centering
        \includegraphics[width=\textwidth]{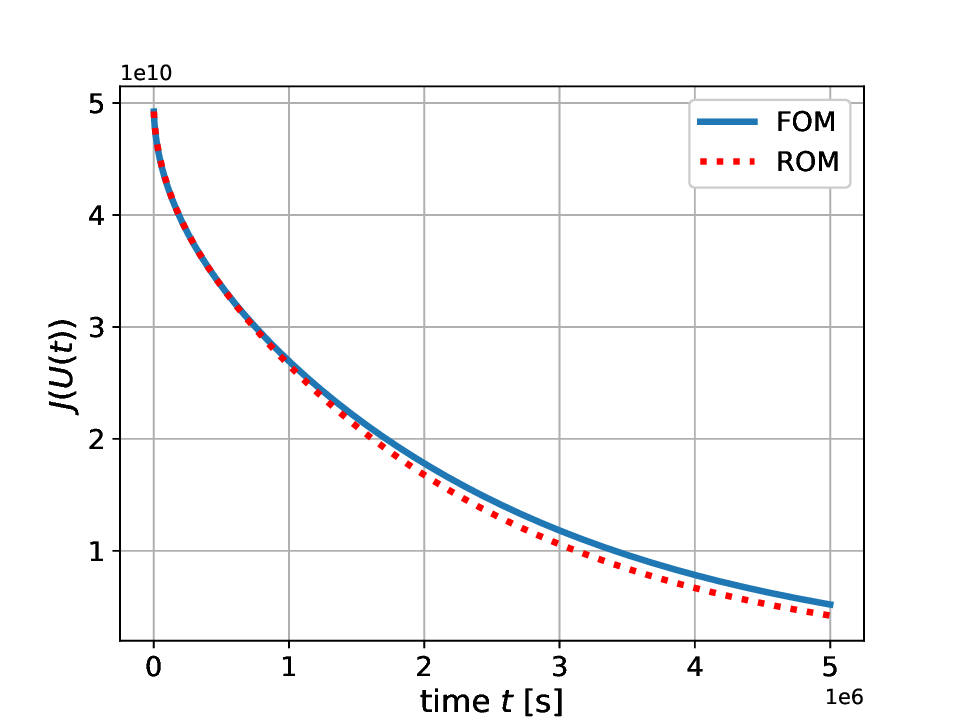}
        \caption{$TOL^{rel} = 5.0\%$}
    \end{subfigure}
    \begin{subfigure}[b]{0.485\textwidth}
        \centering
        \includegraphics[width=\textwidth]{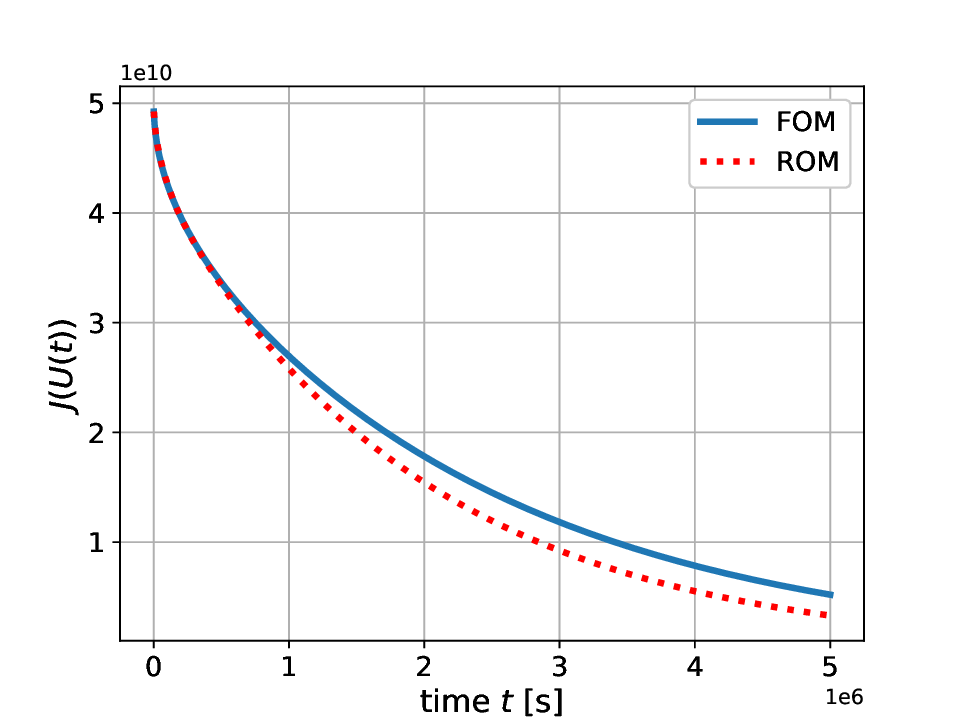}
        \caption{$TOL^{rel} = 10.0\%$}
    \end{subfigure}
    \begin{subfigure}[b]{0.485\textwidth}
        \centering
        \includegraphics[width=\textwidth]{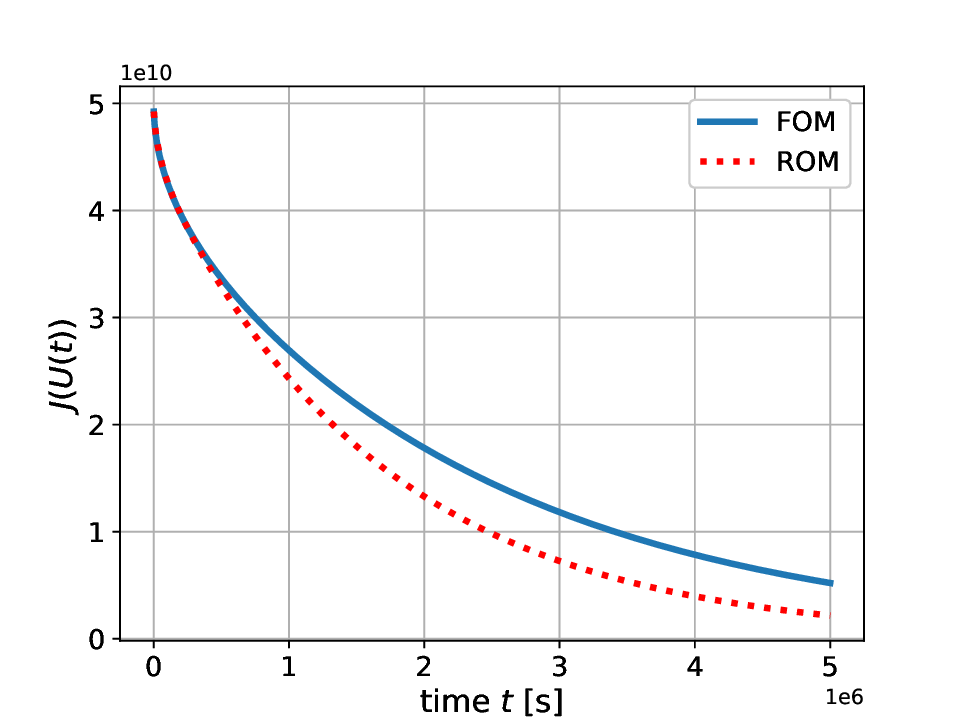}
        \caption{$TOL^{rel} = 20.0\%$}
    \end{subfigure}
    \caption{Goal functionals for Mandel's problem and different relative error tolerances.}
    \label{fig:cost_functional_mandel}
\end{figure}

The time trajectories of the goal functional from the full-order solution $U_h$ and the reduced-order solution $U_N$ are compared for each temporal element in Figure~\ref{fig:cost_functional_mandel}. For these comparisons, we choose relative error tolerances between $0.1 \%$ and $20 \%$.
Both quantity of interest trajectories are almost indistinguishable for tolerances between $0.1 \%$ and $2 \%$. Some differences appear between the full-order and the reduced-order results for larger relative error tolerances. They might be seen as a too large error. However, for certain real-world applications an error of e.g., $10 \%$ might still be acceptable for an appealing computational cost. We can see that choosing the relative error tolerance, the MORe DWR approach allows to control the quality of the approximation of the true quantity of interest, even within a reduced-order solution framework.

Table~\ref{tab:comparison:mandel} gives an overview of the results obtained for the Mandel problem using the MORe DWR algorithm with various error tolerances $TOL^{rel}$  between $0.1\%$ and $20\%$. For each tolerance value, the relative true error $e^{\text{rel}}$ is given as well as the computational speedup, the total  number of FOM solves, the POD basis sizes for the primal displacement/pressure and the dual displacement/pressure, the effectivity index $\Ieff$ from (\ref{eq:effectivity_index}) and the indicator index $\Iind$ defined in (\ref{eq:indicator_index}).
Here, the number of FOM solves sums up all primal and dual solves and the basis sizes are shown in the pattern: primal displacement / primal pressure + dual displacement / dual pressure. 

The relative error as well as the speedup increase with a rise in tolerance. The relative error is rather close to the error tolerance, but sometimes slightly exceeds it. Indeed, here error estimates are employed, not error bounds. This can be explained by the fact that we do not enforce the error tolerance on each temporal element, but only on the whole integrated time domain. In contrast, enforcing the error tolerance on each temporal element leads to much lower error than the error tolerance \cite{FiRoWiChaFau2023}.
The speedup is explained by the decreasing amount of FOM solves and smaller POD bases for both the primal and dual problems, as well as displacement and pressure w.r.t. the given tolerance.  The effectivity index and the indicator index are close to $1$ as expected for a linear PDE and goal functional. This validates the usage of dual-weighted residual error estimates to control the reduced-order-modeling error because the estimated error accurately approximates the true error.

\begin{table}[H]
    \centering
        \begin{tabular}{ |c||c|c|c|c|c|c|c|  }
            \hline
            $TOL^{rel}$ [\%]    &   $e^{rel}$ [\%] &   speedup &   FOM solves & ROM size       &   $\Ieff$ &   $\Iind$ \\
            \hline
  \phantom{0}0.1 & \phantom{0}0.123 &   \phantom{0}8.5 &           97 & 5 / 33 + 37 / 26 &       1.399 &     1.377 \\
  \phantom{0}1\phantom{.0} & \phantom{0}0.821 &   \phantom{0}9.5 &           83 & 5 / 26 + 34 / 25 &       1.089 &     1.092 \\
  \phantom{0}2\phantom{.0} & \phantom{0}2.03\phantom{0}  &  10.1  &           77 & 5 / 24 + 33 / 25 &       1.038 &     1.042 \\
  \phantom{0}5\phantom{.0} & \phantom{0}4.89\phantom{0}  &  11.2  &           69 & 5 / 21 + 31 / 23 &       1.052 &     1.054 \\
  10\phantom{.0} & 10.7\phantom{00}   &  12.7  &           59 & 4 / 16 + 27 / 20 &       1.113 &     1.110 \\
  20\phantom{.0} & 19.6\phantom{00}   &  14.3  &           51 & 4 / 13 + 24 / 18 &       1.143 &     1.142 \\
            \hline
        \end{tabular}
    \caption{Performance of MORe DWR method for the Mandel problem depending on the tolerance in the goal functional.}
    \label{tab:comparison:mandel}
\end{table}

\subsection{Footing problem in 3D}
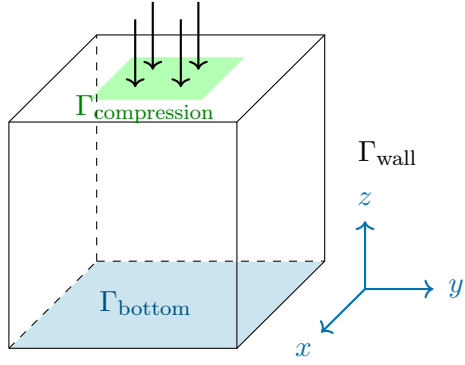
\begin{figure}[H]
    \begin{center}
    \begin{tikzpicture}[scale=3, draw=black]
        \draw[draw=white,fill=blue!20!white] (0,0,0) -- (1,0,0) -- (1,0,-1) -- (0,0,-1) -- (0,0,0);
        \draw[draw=white,fill=green!30!white] (0.25,1,-0.25) -- (0.25,1,-0.75) -- (0.75,1,-0.75) -- (0.75,1,-0.25) -- (0.25,1,-0.25);
        \draw[draw=black] (0,0,0) -- (1,0,0) -- (1,1,0) -- (0,1,0) -- (0,0,0);
        \draw[draw=black] (1,0,0) -- (1,0,-1) -- (1,1,-1) -- (0,1,-1) -- (0,1,0);
        \draw[draw=black] (1,1,0) -- (1,1,-1);
        \draw[draw=black, dashed] (0,0,0) -- (0,0,-1);
        \draw[draw=black, dashed] (0,0,-1) -- (1,0,-1);
        \draw[draw=black, dashed] (0,0,-1) -- (0,1,-1);
        \node (bottom) at (0.6,0.2,0) {\color{blue!70!black}$\Gamma_{\text{bottom}}$};
        \node (compression) at (0.29,0.75,-0.8) {\color{green!50!black}$\Gamma_{\text{compression}}$};
        \node (wall) at (1.2,0.4,-1.2) {$\Gamma_{\text{wall}}$};
        \draw[draw=black, ->, thick]
        (0.4, 1.3, -0.4) -- (0.4, 1.0, -0.4);
        \draw[draw=black, ->, thick]
        (0.4, 1.3, -0.6) -- (0.4, 1.0, -0.6);
        \draw[draw=black, ->, thick]
        (0.6, 1.3, -0.4) -- (0.6, 1.0, -0.4);
        \draw[draw=black, ->, thick]
        (0.6, 1.3, -0.6) -- (0.6, 1.0, -0.6);
        \draw[draw=blue, ->, thick]
        (1.1,-0.2,-1.2) -- (1.4,-0.2,-1.2);
        \node(x) at (1.5, -0.2, -1.2) {\color{blue}$y$};
        \draw[draw=blue, ->, thick]
        (1.1,-0.2,-1.2) -- (1.1,0.1,-1.2);
        \node(z) at (1.1, 0.2, -1.2) {\color{blue}$z$};
        \draw[draw=blue, ->, thick]
        (1.1,-0.2,-1.2) -- (1.1,-0.2,-0.7);
        \node(y) at (1.1, -0.2, -0.5) {\color{blue}$x$};
    \end{tikzpicture}
    \caption{Domain for 3D footing problem.}
    \label{fig:footing_domain}
    \end{center}
\end{figure}

In this second numerical example, a three-dimensional footing problem inspired by \cite{Gaspar2008} is studied. Let $\Omega := (\SI{-32}{\meter}, \SI{32}{\meter}) \times (\SI{-32}{\meter}, \SI{32}{\meter}) \times (\SI{0}{\meter}, \SI{64}{\meter})$ and $I := (\SI{0}{\second}, \SI{5000000}{\second})$ with boundaries as shown in Figure~\ref{fig:footing_domain}.

The initial and boundary conditions are given by
\begin{align*}
    p(0) &= p^0 = 0 &&\text{in } \Omega \times \{0\}, \\
    u(0) &= u^0 = 0 &&\text{in } \Omega \times \{0\}, \\
     \frac{K}{\nu}\nabla_x p\cdot n &= 0 &&\text{on } \partial \Omega \setminus \Gamma_{\text{bottom}} \times I, \quad\tag*{(No flow condition, homogeneous Neumann)}\\
    \sigma(u) \cdot n &= -\bar{t}e_z &&\text{on } \Gamma_{\text{compression}} \times I,\quad\tag*{(Traction condition, inhomogeneous Neumann)} \\
    \sigma(u) \cdot n &= 0 &&\text{on } \Gamma_{\text{top}} \setminus \Gamma_{\text{compression}}\times I,\quad\tag*{(Traction-free, homogeneous Neumann)}\\ 
    p &= 0  &&\text{on } \Gamma_{\text{bottom}} \times I,\quad\tag*{(Constant zero pressure, homogeneous Dirichlet)} \\
    u &= 0  &&\text{on } \Gamma_{\text{bottom}} \times I,\quad\tag*{(Fixed displacements, homogeneous Dirichlet)} \\
    \sigma(u) \cdot n &= 0 &&\text{on } \Gamma_{\text{wall}} \times I, \quad\tag*{(Traction-free, homogeneous Neumann)}
\end{align*}
The material parameters as listed in Table~\ref{tab:params_mandel} are used again.

\indent As before, Taylor-Hood elements, namely quadratic finite elements in space for the displacement $u$ and linear finite elements for the pressure $p$ are employed. 
The spatial mesh is fixed, it comprises 16 spatial cells in each direction i.e., an isotropic mesh with $107,811$ DoFs for displacement and $4,913$ DoFs for pressure.
The temporal domain from $0$ to the end time $T = \SI{5000000}{\second} \approx \SI{58}{\day} $ is discretized with $5,000$ temporal elements, which leads to a 
time step size $k = \SI{1000}{\second}$. For the 
goal functional, we choose the time-integrated pressure acting at the compression boundary i.e.,
\begin{align*}
    J(U) := \int_I \int_{\Gamma_{\text{compression}}} p\ \mathrm{d}x\ \mathrm{d}t.
\end{align*}

A linear system of equations with $112,724$ unknowns is obtained, which has to be solved for each temporal element. The previous 2D Mandel problem was solved with a direct solver. Due to the much larger size of the equation system, we now have to resort to an iterative solution scheme for maintaining a low computational cost. The equation system is not symmetric and the memory for the calculations is not limited; thus we choose the generalized minimal residual method (GMRES). 
Specifically, SciPy \cite{SciPy2020} implementation of GMRES is utilized with a tolerance for convergence of $5 \cdot 10^{-8}$. 
For the preconditioning of the system, we have tested a diagonal Jacobi preconditioner as well as the smoothed aggregation algebraic multigrid (SA-AMG) method utilizing the implementation provided by \cite{pyamg2023}. Both approaches have been benchmarked on the footing problem for the first $500$ temporal elements. The SA-AMG ansatz needs a mean iteration count of $6.33$ for convergence per temporal element solve, while the Jacobi preconditioner requires $283.81$ iterations to reach the same accuracy. Although the number of iterations is significantly lower for the SA-AMG method, the mean wall time per temporal element solve is $93.37\, \mathrm{s}$ for the SA-AMG approach, but only $3.89\, \mathrm{s}$ for the Jacobi preconditioner. The Jacobi preconditioner is significantly faster than the SA-AMG method in our case, most likely due to easier parallelization capabilities on multicore machines. Thus, the Jacobi preconditioner is employed for this application case.

Similar to Mandel's problem, we choose $1-10^{-7}$ and $1-10^{-11}$ for the primal displacement and pressure energy tolerance \eqref{eq:energy_content}.
Again, the dual energy tolerances are both set to $1-10^{-9}$. As argued in Section~\ref{sec:incremental_ROM}, the dual spaces are enriched with additional dual full-order solutions. In this case, we perform this additional enrichment for the first 8 MORe DWR iterations and include the dual solutions for the last five time steps of the dual problem i.e., the time steps $[t_0, t_4]$ of the primal problem discretization.

\begin{figure}[H]
    \centering
    \begin{subfigure}[b]{0.485\textwidth}
        \centering
        \includegraphics[width=\textwidth]{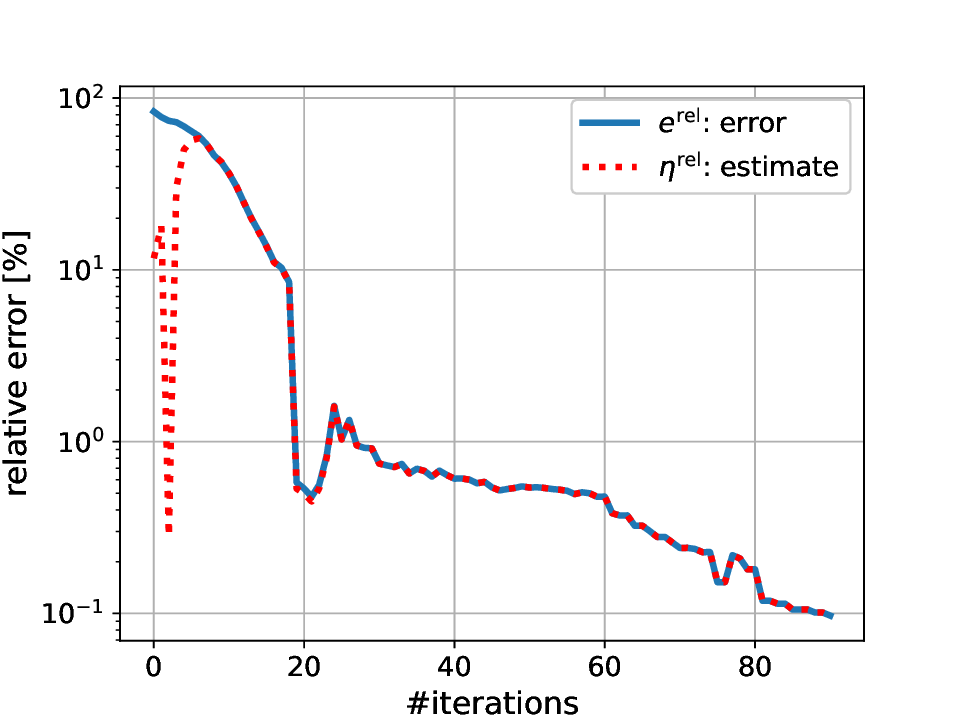}
        \caption{Evolution of the true and estimated relative errors throughout the MORe DWR iterations.}
        \label{fig:adaptivity_footing:iteration_error}
    \end{subfigure}
    \hfill
    \begin{subfigure}[b]{0.485\textwidth}
        \centering
        \includegraphics[width=\textwidth]{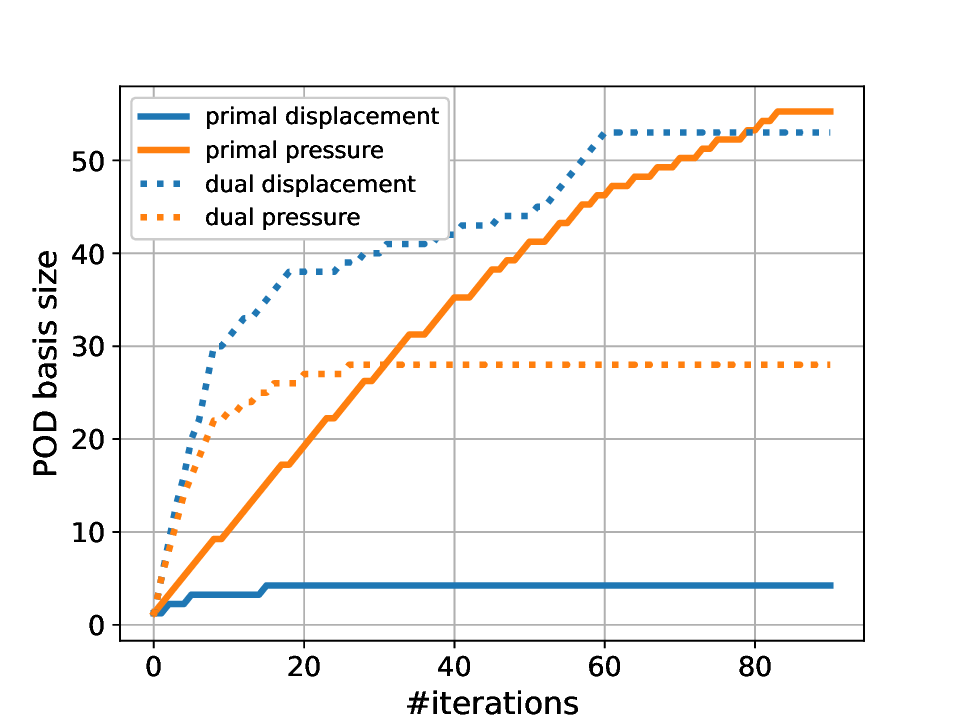}
        \caption{Evolution of the POD basis sizes throughout the MORe DWR iterations. \phantom{this is a placeholder}}
        \label{fig:adaptivity_footing:pod_basis}
    \end{subfigure}
    \caption{Adaptivity in the MORe DWR method of the 3D footing problem.}
    \label{fig:adaptivity_footing}
\end{figure}

The adaptive enrichment of the method is illustrated in Figure~\ref{fig:adaptivity_footing}. The relative error between the full-order and reduced-order solutions, as well as the error estimator is shown in Figure~\ref{fig:adaptivity_footing:iteration_error}. 
We observe a similar behavior for the first iterations as in Mandel's problem: a discrepancy between the estimate and the error, with a severe underestimation of the error. Thereafter, both quantities align for the rest of the iterations. In addition, we see a sharp decline after 18 iterations that is followed by a non-monotonic behavior, highlighting the greedy enrichment process of the MORe DWR method.

The evolution of the size of the primal/dual POD bases for displacement/pressure is shown in Figure~\ref{fig:adaptivity_footing:pod_basis}. We see again a steep increase in the dual POD sizes in the first iterations that then flattens out and is finally surpassed by the primal pressure POD size. The size of the primal displacement POD size is almost constant, with a maximum of 4.

\begin{figure}[H]
    \centering
    \begin{subfigure}[b]{0.485\textwidth}
        \centering
        \includegraphics[width=\textwidth]{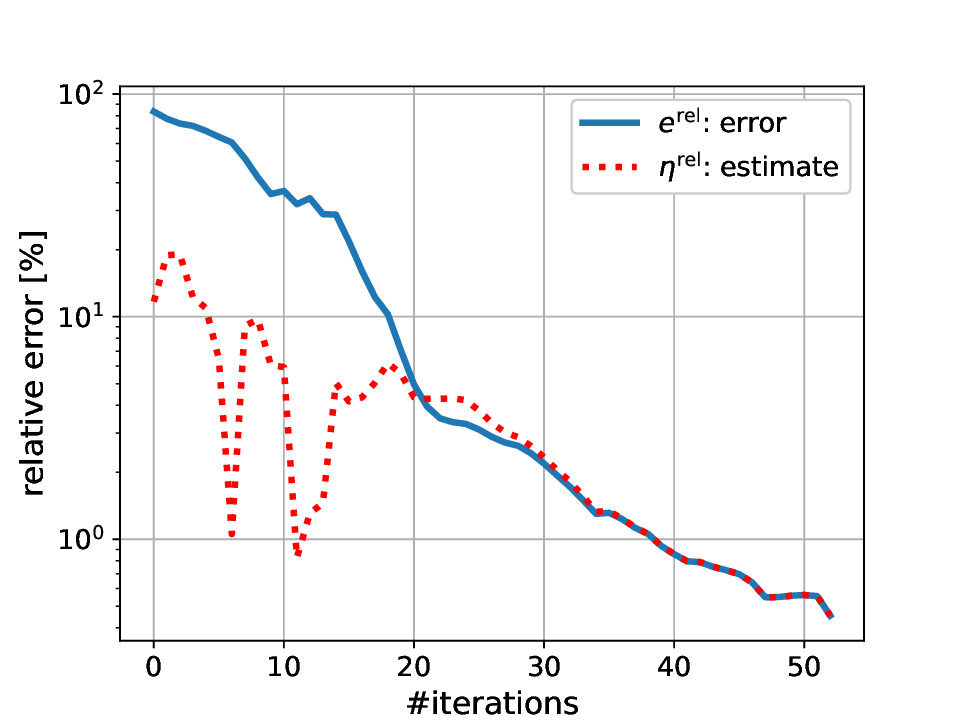}
        \caption{Evolution of the true and estimated relative errors throughout the MORe DWR iterations.}
        \label{fig:adaptivity_footing_without_dual:iteration_error}
    \end{subfigure}
    \hfill
    \begin{subfigure}[b]{0.485\textwidth}
        \centering
        \includegraphics[width=\textwidth]{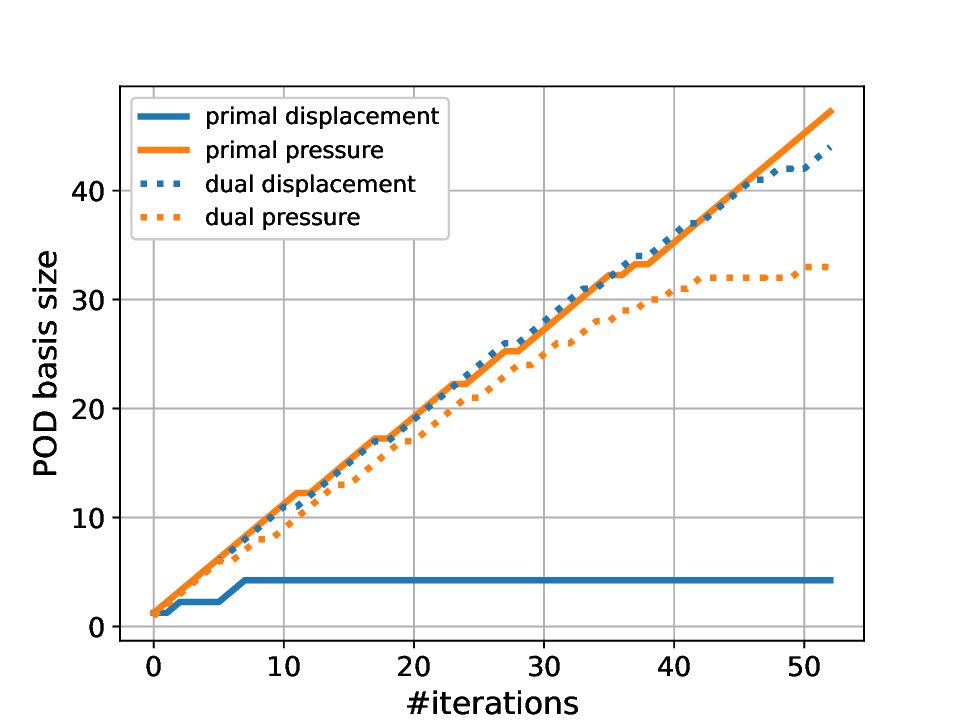}
        \caption{Evolution of the POD basis sizes throughout the MORe DWR iterations. \phantom{this is a placeholder}}
        \label{fig:adaptivity_footing_without_dual:basis_size}
    \end{subfigure}
    \caption{Adaptivity in the MORe DWR method of the 3D footing problem \textit{without additional dual basis enrichments}.}
    \label{fig:adaptivity_footing_without_dual}
\end{figure}

Next, in Figure~\ref{fig:adaptivity_footing_without_dual} we compare the latter results with the MORe DWR method \textit{neglecting} the additional dual reduced space enrichment proposed in Section~\ref{sec:incremental_ROM}. A minimum number of $20$ iterations is enforced to ensure the consistency of the results because it is expected from the previous study that the error is severely underestimated during the first iterations without extra enrichment (Figure~\ref{fig:adaptivity_footing:iteration_error}).
Subfigure~\ref{fig:adaptivity_footing_without_dual:iteration_error} compares the relative error between the full-order and reduced-order solutions with its estimate over the course of MORe DWR iterations. In contrast to the results in Figure~\ref{fig:adaptivity_footing:iteration_error}, we observe that the first phase of discrepancy between the estimate and the error, with the typical severe underestimation of the error, lasts much longer with more variation. The true error and its estimate align after 35 iterations. Thereafter, the error estimate is as accurate as when using the extra dual enrichment. Thus, the extra enrichment of the dual basis appears necessary for an accurate error estimate from the early stage of the algorithm. 
The ROM solution that complies to a tolerance of $0.5\%$ is obtained with a total of $104$ FOM evaluations (52 primal, 52 dual evaluations) without the dual enrichment, whereas only 80 FOM evaluations (20 primal, 20 dual, 40 extra dual evaluations) were required with the enriched dual basis.
Thus, although we added dual FOM evaluations in the first iterations of the extra-enriched approach, the trade-off seems worth it, as this approach reduces the total number of iterations and thus the computational costs.

In addition, the development of the POD bases in case of no additional dual enrichment is 
shown in Figure \ref{fig:adaptivity_footing_without_dual:basis_size}.
In contrast to Figure~\ref{fig:adaptivity_footing:pod_basis}, except for the primal displacement, all the POD bases start growing at nearly the same rate. The dual pressure basis size stagnates after 42 iterations, and the primal displacement has, after 7 iterations, a constant size of 4. Comparing the POD basis sizes for both approaches  for a relative tolerance of $0.5\%$, the final dual POD basis sizes are nearly the same. In contrast, in the case of the extra enrichment, the primal pressure POD size of 19 is considerably smaller than the 47 modes needed when neglecting this enrichment.

\begin{figure}[H]
    \centering
    \begin{subfigure}[b]{0.485\textwidth}
        \centering
        \includegraphics[width=\textwidth]{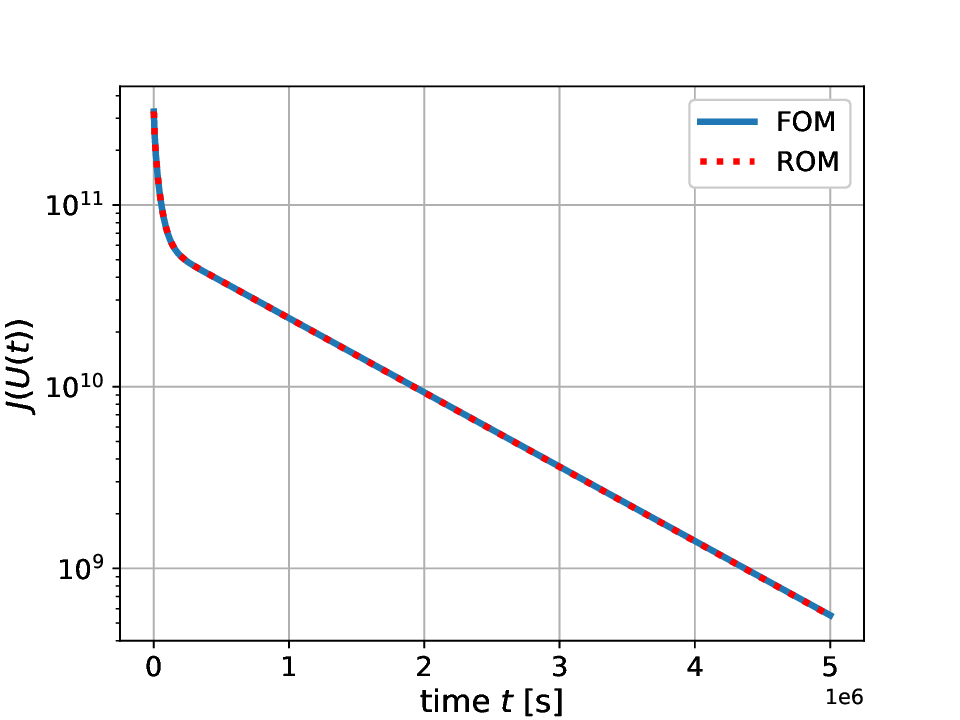}
        \caption{$TOL^{rel} = 0.1\%$}
    \end{subfigure}
    \hfill
    \begin{subfigure}[b]{0.485\textwidth}
        \centering
        \includegraphics[width=\textwidth]{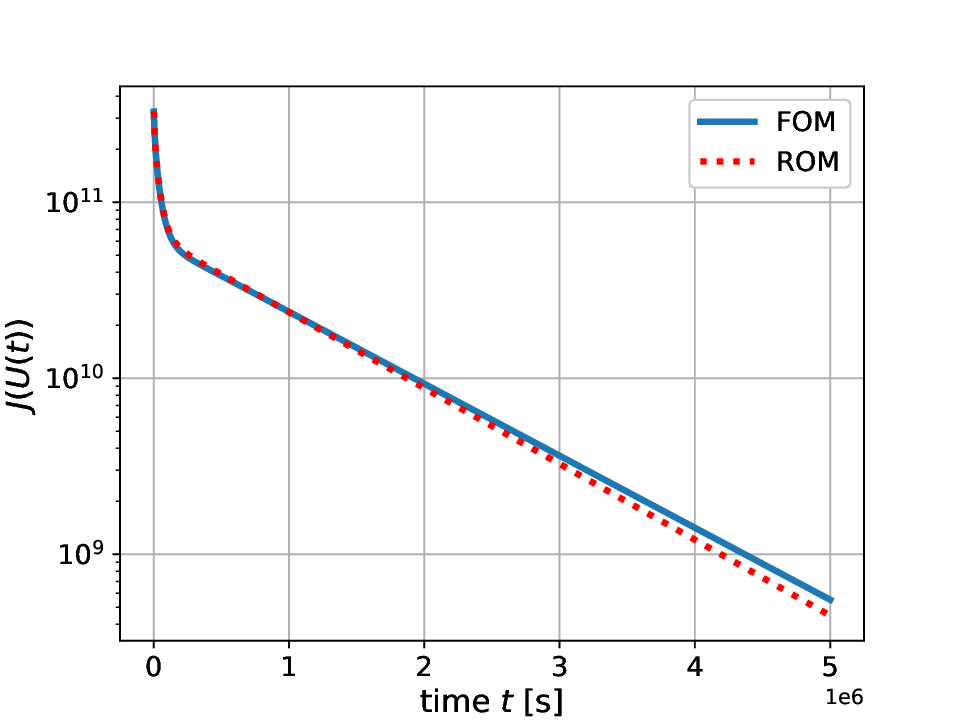}
        \caption{$TOL^{rel} = 1.0\%$}
    \end{subfigure}
    \begin{subfigure}[b]{0.485\textwidth}
        \centering
        \includegraphics[width=\textwidth]{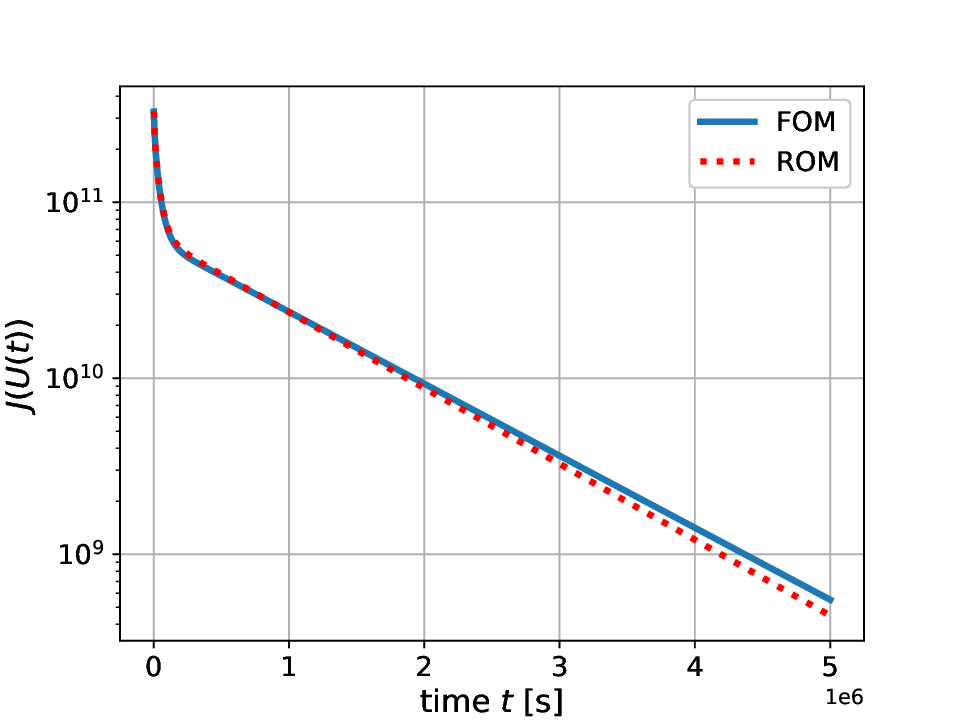}
        \caption{$TOL^{rel} = 2.0\%$}
    \end{subfigure}
    \begin{subfigure}[b]{0.485\textwidth}
        \centering
        \includegraphics[width=\textwidth]{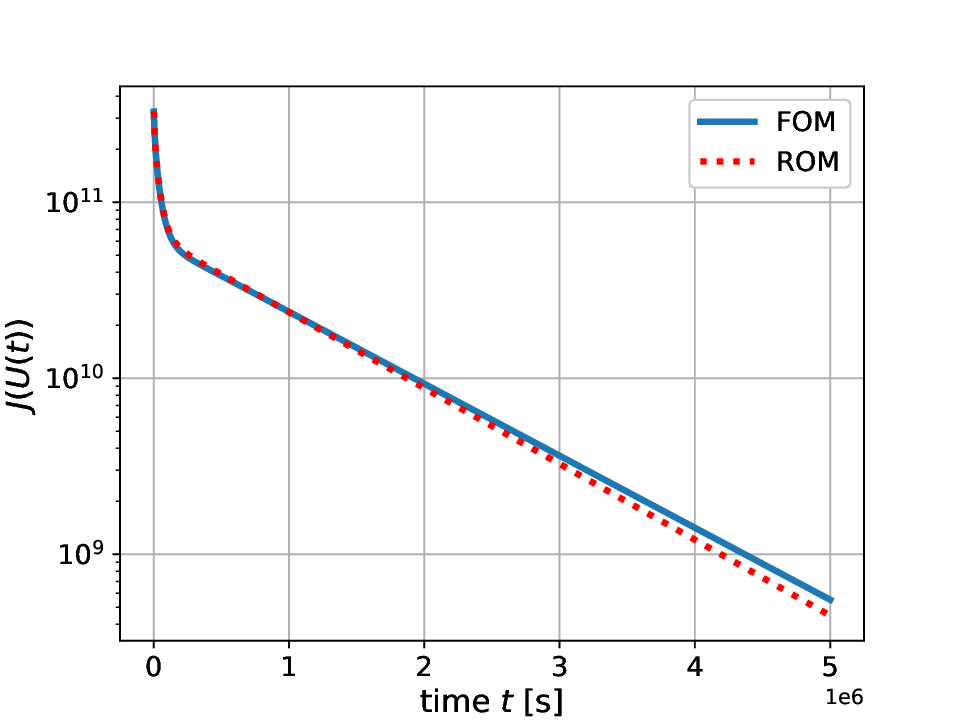}
        \caption{$TOL^{rel} = 5.0\%$}
    \end{subfigure}
    \begin{subfigure}[b]{0.485\textwidth}
        \centering
        \includegraphics[width=\textwidth]{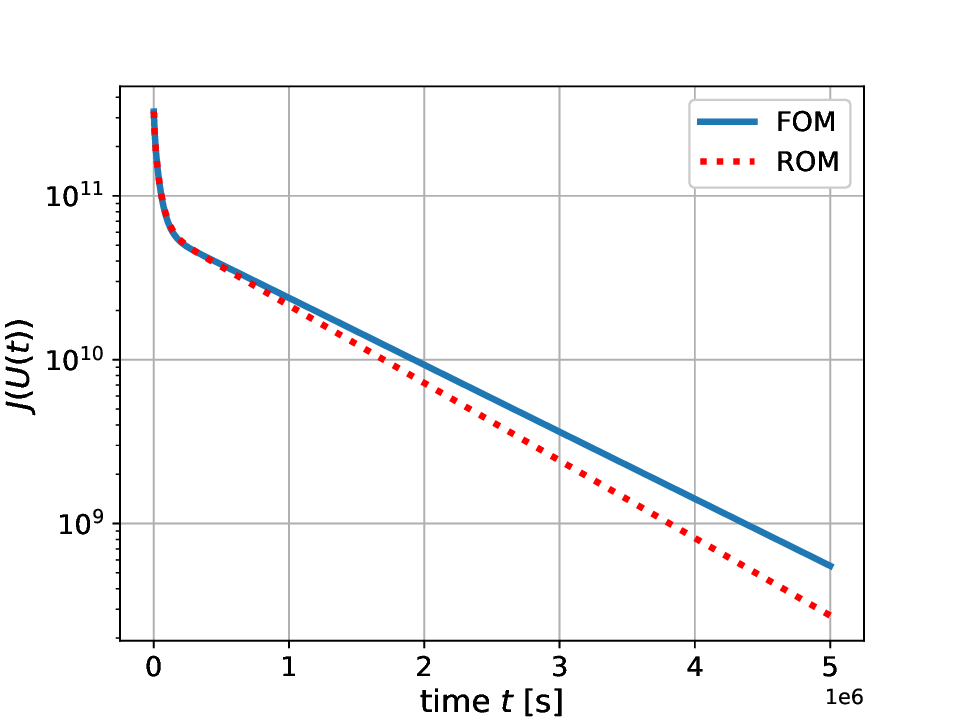}
        \caption{$TOL^{rel} = 10.0\%$}
    \end{subfigure}
    \begin{subfigure}[b]{0.485\textwidth}
        \centering
        \includegraphics[width=\textwidth]{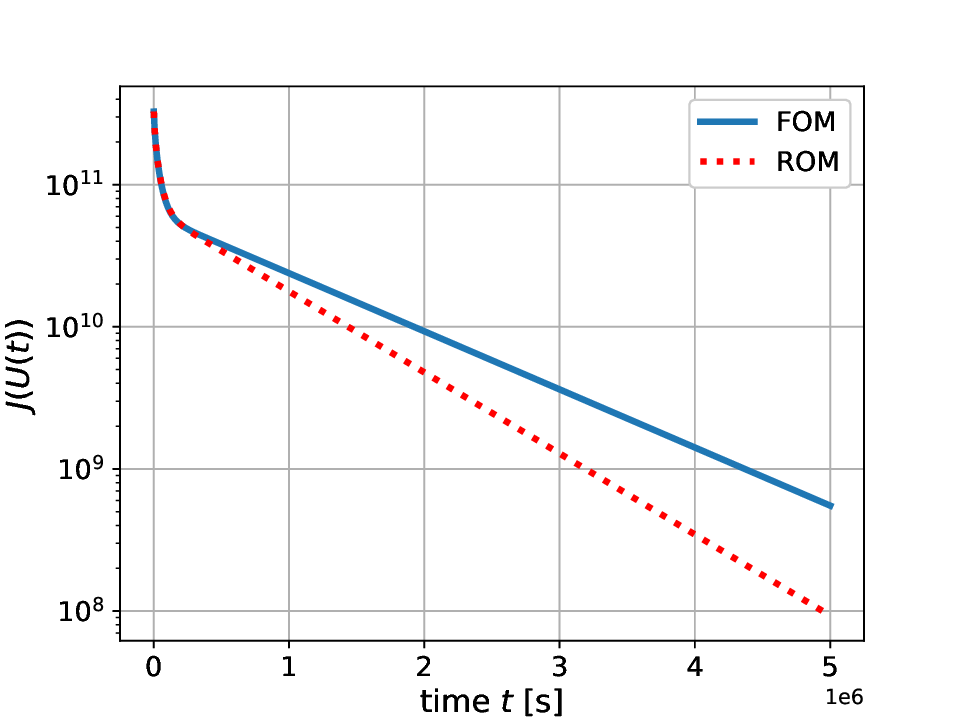}
        \caption{$TOL^{rel} = 20.0\%$}
    \end{subfigure}    
    \caption{Goal functionals for the 3D footing problem, depending on the tolerance in the goal functional. 
    }
    \label{fig:cost_functional_footing}
\end{figure}

The time trajectories of the goal functional from the reduced-order solution is compared with the trajectory given by the full-order solution for error tolerances  between 0.1\% and 20\% in Figure~\ref{fig:cost_functional_footing}. The trajectories of  quantities of interest are nearly identical for a tolerance of 0.1\%. Differences emerge between the full-order and reduced-order goal functionals for larger tolerances. In essence, the reduced-order solution reaches the accuracy chosen by the user.

The time trajectories of the goal functional from the reduced-order solution is compared with the trajectory given by the full-order solution for error tolerances  between 0.1\% and 20\% in Figure~\ref{fig:cost_functional_footing}. The trajectories of both quantities of interest are nearly identical for a tolerance of 0.1\%. Differences emerge between the full-order and reduced-order goal functionals for larger tolerances. In essence, the reduced-order solution reaches the accuracy chosen by the user.

In Table~\ref{tab:comparison:footing}, we give an overview of the simulation results of the three-dimensional footing problem for different error tolerances  between $0.1\%$ and $20\%$.
The tolerances are met besides in the case of $0.5\%$, where the error is slightly underestimated. We observe that with a rise in the tolerance, the relative error as well as the speedup increase. We note that the algorithm has the same behavior for the tolerances $1\%$, $2\%$ and $5\%$. Indeed, these tolerances are reached within the same number of MORe DWR iterations, yielding all other quantities to be the same. This behavior can be explained by the sharp decline in the error after 19 MORe DWR iterations, c.f. Figure~\ref{fig:adaptivity_footing:iteration_error}. The evaluation of the dual-weighted residual error estimates by means of the effectivity and indicator indices shows that the estimate integrated over the temporal domain is almost perfect, as the effectivity index is near $1.0$ in all cases. However, the indicator indices vary from $1.00$ up to $3.44$ which translates to inaccurate temporal error localizations in the case of high indicator indices.

\begin{table}[H]
    \centering
        \begin{tabular}{|c||c|c|c|c|c|c|c|}
\hline
        $TOL^{rel}$ [\%]    &   $e^{rel}$ [\%] &   speedup &   FOM solves & ROM size       &   $\Ieff$ &   $\Iind$ \\
\hline
  \phantom{0}0.1 &             \phantom{0}0.0971 &   \phantom{2}8.6 &           220& 4 / 55 + 53 / 28 &       0.999   &     1.207 \\
  \phantom{0}0.5 &             \phantom{0}0.5333 &   21.2 &           80 & 4 / 19 + 38 / 27 &       1.068   &     3.441 \\
  \phantom{0}1\phantom{.0}   &           \phantom{0}0.579\phantom{0} &   22.4 &           78 & 4 / 18 + 38 / 26 &       1.084 &     3.378 \\
  \phantom{0}2\phantom{.0}   &           \phantom{0}0.579\phantom{0} &   21.7  &           78 & 4 / 18 + 38 / 26 &       1.084 &     3.378 \\
  \phantom{0}5\phantom{.0}   &           \phantom{0}0.579\phantom{0} &   22.2 &           78 & 4 / 18 + 38 / 26 &       1.084 &     3.378 \\
 10\phantom{.0}   &             \phantom{0}8.49\phantom{00}   &   22.4 &           76 & 4 / 17 + 38 / 26 &       1.008 &     1.099   \\
 20\phantom{.0}   &             19.9\phantom{000}   &   26.2 &           66 & 3 / 13 + 33 / 24 &       1.005 &     1.031 \\
\hline
\end{tabular}
    \caption{Performance of MORe DWR method for the 3D footing problem, depending on the tolerance in the goal functional.}
    \label{tab:comparison:footing}
\end{table}

\section{Conclusion}

We developed an adaptive incremental reduced-order (MORe DWR) model for single-phase flow problems in porous media, namely the Biot system. To this end, we further developed and extended 
our prior work \cite{FiRoWiChaFau2023}. Therein, the following ingredients are combined. 
First, a space-time formulation and Galerkin space-time finite element discretization 
for the coupled problem. Second, a ROM method based on an incremental POD. Third, 
a goal-oriented a posteriori error estimator using the dual-weighted residual method for 
adaptively refining the ROM bases.
The efficiency of the methodology has been demonstrated, by using the incremental ROM for the two-dimensional Mandel benchmark and a three-dimensional footing problem. The algorithm is stopped according to the total relative error considering the whole temporal domain. The speed-ups offered by this approach compared with the reference full-order approach vary from 8.5x in two dimensions for the smaller relative tolerances, to 26.2x for the three-dimensional case with large relative tolerances. 
The accuracy of the error estimator has been improved by additional enrichment of the dual reduced bases. The true error and the error estimator almost coincide with effectivity indices close to one, which makes the error estimates reliable. It is important to note that the reduced-order model starts with no other prior information than a single primal snapshot and a dual FOM snapshot.
Hence, the full-order model needs to be queried in the online stage to gather snapshots for basis enrichment. Consequently, our approach could be further sped up by efficient full-order model solvers.
Therefore, for the three-dimensional test, we used GMRES as our iterative solver and compared a simple diagonal Jacobi preconditioner and an algebraic multigrid preconditioner.  
An aspect for future work is the combination of temporal and reduced basis adaptivity. The MORe DWR algorithm already requires the solution of the primal and dual problems. This additional computational effort can be reused to refine the temporal mesh adaptively using the dual-weighted residual method and could reduce the number of time steps, thus further speeding up the computations.

\section{Acknowledgements}
The authors acknowledge the funding of the German Research Foundation (DFG) within the framework of the International Research Training Group on  Computational Mechanics Techniques in High Dimensions GRK 2657 under Grant Number 433082294. The support of the French-German University through the French-German Doctoral college "Sophisticated Numerical and Testing Approaches" (CDFA-DFDK 19-04) is also acknowledged. Moreover, T. Wick gratefully acknowledges the J. Tinsley Oden Faculty Fellowship Research Program at UT Austin for the visit in July/August 2023.


\section{Credit Authorship Contribution Statement}
\makebox[2.75cm][l]{\textbf{H. Fischer}:} Conceptualization, Methodology, Software, Validation, Formal analysis, \\ 
\makebox[2.75cm][l]{} Investigation, Writing – original draft, Writing – review \& editing, Visualization. \\
\makebox[2.75cm][l]{\textbf{J. Roth}:} Conceptualization, Methodology, Software, Validation, Formal analysis,
\\ 
\makebox[2.75cm][l]{} Investigation, Writing – original draft, Writing – review \& editing, Visualization. \\
\makebox[2.75cm][l]{\textbf{L. Chamoin}:} Conceptualization, Formal analysis, Writing - review \& editing, Supervision, Funding 
\\ 
\makebox[2.75cm][l]{} acquisition. \\
\makebox[2.75cm][l]{\textbf{A. Fau}:} Conceptualization, Formal analysis, Writing - review \& editing, Supervision, Funding 
\\ 
\makebox[2.75cm][l]{} acquisition. 
\\
\makebox[2.75cm][l]{\textbf{M. Wheeler}:} Conceptualization, Formal analysis, Discussions, Writing - review \& editing.\\
\makebox[2.75cm][l]{\textbf{T. Wick}:} Conceptualization, Formal analysis, Resources, Writing - review \& editing,
\\ 
\makebox[2.75cm][l]{} Supervision, Funding acquisition. \\



\end{document}